%% file: main.tex
\magnification=1200
\parindent=0pt

\input psfig

%\font\baton msym10
%\def\Z{{\hbox{\baton Z}}}
%\def\N{{\hbox{\baton N}}}
%\def\Q{{\hbox{\baton Q}}}
%\def\R{{\hbox{\baton R}}}

\def\N{{\bf N}}
\def\Q{{\bf Q}}

\def\Z{{\bf Z}}
\def\C{{\bf C}}
\def\T{{\bf T}}

\def\semi{\;;\;}
\def\colon{\;:\;}

\def\A{{\cal A}}
\def\B{{\cal B}}
\def\CC{{\cal C}}
\def\D{{\cal D}}

\def\L{{\cal L}}
\def\barO{\overline {\cal O}}
\def\P{{\cal P}}
\def\RR{{\cal R}}
\def\QQ{{\cal Q}}
\def\card{\mathop{\rm Card}}

\def\xs{X_\sigma} 
\def\ts{T_\sigma}
\def\xtau{X_\tau} 
\def\ttau{T_\tau}
\def\xb{X_\B}
\def\vb{V_\B}
\def\xz{X_\zeta}

\def\card{\mathop{\rm Card}}

% ----------------------------------------------
% COMPILATION CONDITIONNELLE
\long\def\STOP#1\GO{{\bf\par Ici passage sautŽ\par}}
%--------------------------------------------------

%       TITRES

\font\petit=cmr8
\font\gros=cmbx12
\newif\ifpagetitre \pagetitretrue

\newtoks\hautpage \hautpage={\hfil}
\newtoks\hautgarde \hautgarde={\hfil}
\newtoks\baspage \baspage={\hfil\tenrm\folio\hfil}
\newtoks\basgarde \basgarde={\hfil \tenrm \folio\hfil}

\headline={\ifpagetitre\the\hautgarde\else\the\hautpage\fi}
\footline={{\ifpagetitre\the\basgarde\else\the\baspage\fi}%
\global\pagetitrefalse}
 %-----------------------------------------------------------

\def\cqfd{\vbox{\hrule%
\hbox{\vrule height 5pt depth 0pt\kern 5pt\vrule height 5pt depth 0pt}%
\hrule}}

\def\xmin{x_{{\hbox{\petit min}}}}

\def\ruv{R_{u.v}}
\def\rruv{\RR_{u.v}}
\def\duv{\D_{u.v}}
\def\phiuv{\phi_{u.v}}

\def\um{{\,\overrightarrow m}}
\def\one{\hbox{\bf 1}}

\def\th{{\hbox{\petit th}}}

\def\fleche #1{\buildrel{#1}\over\longrightarrow}

% ----------------------------
\centerline{\gros SUBSTITUTION DYNAMICAL SYSTEMS,}
\centerline{\gros BRATTELI DIAGRAMS}
\centerline{\gros AND DIMENSION GROUPS.}

\vskip 0,5cm

\centerline{{\bf F. Durand}\footnote{$^\dagger$}{University of Aix--Marseille II and Institut de Math\'ematiques de Luminy--CNRS, France.}, {\bf B. Host}\footnote{$^\ddagger$}{University of Marne la Vall\'ee, France.}, {\bf C. Skau}\footnote{*}{Norwegian University of Science and Technology, Trondheim, Norway.}}

\vskip 0,5cm

\bigskip
\input section0.tex

\bigskip
\input section1.tex

\bigskip
\input section2.tex

\bigskip
\input section3.tex

\bigskip
\input section4.tex

\bigskip
\input section5.tex

\bigskip
\input section6.tex

\bigskip
\input section7.tex

\bigskip
\input section8.tex

\bigskip
\input biblio.tex

\bye

%% file: section0.tex
{\gros 0. Introduction.}

\bigskip

The present paper explores substitution minimal systems and their relation
to stationary Bratteli diagrams and stationary dimension groups. The
constructions involved are algorithmic and explicit, and render an
effective method to compute an invariant of (ordered) $K$-theoretic nature
for these systems. This new invariant is independent of spectral invariants
which have previously been extensively studied. Before we state the main 
results we give some background.

\medskip

{\bf From Bratteli diagrams to topological dynamics.} In 1972 O. Bratteli [Br] introduced special infinite graphs --- subsequently
called {\sl Bratteli diagrams} --- which conveniently encoded the successive
embeddings of an ascending sequence $(A_n ; n\geq 0)$ of finite-dimensional semi-simple algebras over $\C$ (``multi-matrix algebras''). The sequence $(A_n;n\in\N)$ determines a so-called approximately finite-dimensional (AF) $C^*$-algebra. O. Bratteli proved that the equivalence relation on Bratteli diagrams generated by the operation of telescoping is a complete isomorphism
invariant for AF-algebras. Inspired by O. Bratteli's results, G. A. Elliott [El] introduced in 1976 the notion of {\sl dimension groups}. In fact, to a Bratteli diagram is naturally associated a sequence of free abelian groups of finite rank (simplicially ordered) with order-preserving homomorphisms between successive groups. The resulting inductive limit --- endowed with the induced ordering --- is the dimension group associated
to the Bratteli diagram. The $K$-theoretic underpinning of this notion was 
soon realized, being in fact order-isomorphic to the (ordered)
$K_0$-group of the AF-algebra determined by the Bratteli diagram. Bratteli's
result could now be given a more succinct formulation: a complete 
isomorphism invariant for AF-algebras is the associated dimension group
(as an ordered group with distinguished order unit). Then, in 1980, E. G. Effros, D. Handelman and C.-L. Shen [EHS] gave an axiomatic characterization of dimension groups, namely as countable, ordered abelian groups which are unperforated and satisfy the Riesz interpolation property. Their result made it possible to easily produce explicit examples of dimension groups and to make certain inroads into the classification problem.

From a different direction came the extremely fruitful idea of A. M. Vershik [Ve] to associate dynamics (``adic transformations'') with Bratteli diagrams
(``Markov compacta'') by introducing a lexicographic ordering on the infinite
paths of the diagram. By a careful refining of Vershik's construction,
R. H. Herman, I. F. Putnam and C. F. Skau [HPS] succeeded in showing that every Cantor minimal
system is isomorphic to a Bratteli-Vershik system, i.e. the dynamical system
associated to a properly ordered Bratteli diagram. An immediate consequence
of their result is that every (simple) dimension group occur as the (ordered)
$K_0$--group that is associated to a Cantor minimal system, the latter notion
being defined in purely dynamical terms. Subsequent work by T. Giordano, I. Putnam and C. F. Skau [GPS] showed that the orbit structure of a Cantor minimal system is intimately related to it's $K_0$-group.

\medskip

The main motivation of this paper is that up to now there has been few explicit and natural examples that have been 
worked out to illustrate the correspondence between Cantor minimal systems
and dimension groups, respectively Bratteli diagrams. In this paper we shall
show that the family of substitution minimal systems will yield examples of
the desired nature.

\bigskip

General notions of substitutions appeared for the first time in 1963 in a paper of W. H. Gottshalk [Go]. He pointed out that the famous Morse sequence, initially introduced to prove the existence of recurrent geodesics on surfaces of negative curvature, can be defined by a substitution, now known as the ``Morse substitution''.

\medskip

Nowadays substitutions arose and are studied in many domains. We briefly recall some of these domains.

\medskip

{\bf Substitutions and topological dynamics.} The study of symbolic dynamical systems has been initiated by G. A. Hedlund and M. Morse in [HM1,HM2]. One of the simplest ways of constructing interesting minimal symbolic dynamical systems is by means of substitutions. For this reason, ergodic and topological properties of substitution systems has been extensively studied (see [Qu] for a good bibliography on this subject). Numerous papers deal with the spectrum of substitution system. First results were obtained for substitutions of constant length [De2,Ma] and then for substitutions of non-constant length [FMN,Ho,So].

\medskip

{\bf D0L systems.} Independently substitutions (called iterated morphisms) appeared to be attractive objects for the theory of formal languages; where they are studied through the notion of D0L systems. G. Rozenberg and A. Salomaa [RS] proved the decidability of the D0L equivalence problem.

\medskip

{\bf $p$-recognizable sets of numbers.} Cobham's theorem [Co1] belongs to the domains of the theory of automata and of Arithmetic. It states that a set of positive integers which is ``recognizable'' in two multiplicatively independent bases is eventually periodic. Sets of integers which are recognizable in base $p$ (i.e. $p$-recognizable) are characterized by means of substitutions of constant length $p$ [Co2]. This gives a formulation of Cobham's theorem using substitutions of constant length. 

Let $p$ be a power of a prime number, in [CKMR] it is shown that a set is $p$-recognizable if and only if the power series associated to its characteristic sequence is algebraic over some field of rational fractions.

These results are connected to some problems of transcendence. For example F. M. Dekking [De1] proved that the number $\sum_{i=1}^{+\infty} x_n 2^{-n}$, where $(x_n ; n\in \N)$ is the Morse sequence, is transcendent. Articles [LP] and [FM] give a partial answer to the following question: What can we say about real numbers having an expansion in some base $p$ which is substitutive?

\medskip

A characterization of sequences generated by substitutions of constant length  were given in [BHMV] with the help of the formalism of first order logic.

\medskip

{\bf Geometry.} The article [Ra] of G. Rauzy opened a new area of investigation for substitutions. He constructed a partition of the torus $\T^2$ in three fractal sets, and a rotation on $\T^2$, and looked at the symbolic system obtained in coding the orbits of this rotation by the partition. He found that this system is isomorphic to some substitution subshift. The partition G. Rauzy introduced generates self-similar tilings. In [De3] F. M. Dekking gave a method, using substitutions, to generate fractal curves and self-similar tilings; in particular it generates the famous Penrose tiling [Pe].

\medskip

{\bf Theoretical physics.} In 1984 D. Shechtman and al. [SBGC] discovered the existence in nature of quasi-periodic tilings: the quasicrystals. Later substitution sequences appeared to be good one-dimensional models for quasicrystals (see [BT]) and were used in a lot of papers of theoretical physics. In particular, several authors have studied the properties of the discrete Schr\"{o}dinger operator with potential given by a substitution sequence (see [Be,BBG]).

\medskip

{\bf Contents.} Before we state the main results of this paper we want to acknowledge the
paper by A. H. Forrest [Fo], where the main part of Theorem 1 is proved. However, the proofs given in his paper are mostly of existential nature and do
not state a feasible method to compute effectively dimension groups and Bratteli diagrams associated to substitution systems (see also [AP]).

\medskip

{\bf Theorem 1.} {\sl The family $\B$ of Bratteli-Vershik systems associated to stationary,
properly ordered Bratteli diagrams is (up to isomorphism) the disjoint union
of the family of substitution minimal systems and the family  of stationary
odometer systems. Furthermore, the correspondence in question is given by
an explicit and algorithmic effective construction. The same is true for
the computation of the (stationary) dimension group associated to a 
substitution minimal system.}

\medskip

By invoking the characterization of Kakutani equivalence of Cantor minimal
systems in terms of ordered Bratteli diagrams of [GPS; Theorem 3.8] we get 
the following corollary.

\medskip

{\bf Corollary 2.} {\sl The family $\B$ is stable under Kakutani equivalence.}

\medskip

In this paper we introduce the notion of linearly recurrent subshifts and we prove the following result.

\medskip

{\bf Theorem 3.} {\sl Let $(Y,T)$ be a linearly recurrent aperiodic subshift. There exists a constant $D$ such that the length of each proper chain of aperiodic subshift factor maps 
$$ 
(Y_{n},T)\; \fleche{ \gamma_{n-1}} \;
(Y_{n-1},T)\;  \fleche{ \gamma_{n-2}}\; 
\cdots\; \fleche{\gamma_1}(Y_1,T)\;
\fleche{\gamma_0} \; (Y_0,T)\  ,
$$
with $(Y_n,T) = (Y,T)$, is less than $D$ (i.e. $n\leq D$).
} 

\medskip

We show that substitution minimal systems are linearly recurrent which is helpful, using Theorem 3 also, to prove the following theorem.

\medskip

{\bf Theorem 4.} {\sl A Cantor factor of a system (X,T) belonging to the family $\B$ in the above theorem again belongs to $\B$.}

\medskip

%% file: section1.tex
{\gros 1. Bratteli diagrams and dimensions groups -- Ordered
Bratteli diagrams and Cantor minimal systems -- Stationary diagrams.}

\bigskip

In this section we will recall some basic concepts and results, and
also some of the key constructions that are used to obtain these
results.  This will later be applied in the context of stationary
Bratteli diagrams and substitution minimal systems --- the focus of
this paper.  We will introduce some relevant terminology and
notation.  We refer to [HPS,GPS,Ef] for further details.

\medskip

{\bf 1.1. Bratteli diagrams.}

\medskip

{\bf 1.1.1. Definition 1.}  {\sl A {\rm Bratteli diagram} is an infinite
directed graph $(V, E)$, such that the vertex set $V$ and the edge
set $E$ can be partitioned into finite sets

$$
V = V_0 \cup V_1 \cup V_2 \cup \cdots \ \ {\rm and} \ \ E = E_1 \cup E_2 
\cup \cdots 
$$

with the following properties:

\medskip

$\imath )$ $V_0 = \{ v_0 \}$ is a one-point set.

$\imath\imath )$ $r(E_n) \subseteq V_n, \; s(E_n) \subseteq V_{n-1}, n = 1,2,
\ldots,$ where $r$ is the associated {\rm range map} and $s$ is the associated 
{\rm source map}.  Also, $s^{-1} (v) \neq \emptyset$  for all $v \in V$ and 
$r^{-1} (v) \neq \emptyset$  for all $v \in V \setminus V_0$.}

\medskip

There is an obvious notion of isomorphism between Bratteli diagrams $(V, E)$
and $(V' , E' )$~; namely, there exist a pair of bijections between $V$
and $V'$ and between $E$ and $E'$, respectively, preserving the gradings and
intertwining the respective source and range maps.

It is convenient to give a diagrammatic presentation of the Bratteli
diagram with $V_n$ the vertices at (horizontal) level $n$, and $E_n$
the edges (downward directed) connecting the vertices at level $n-1$
with those at level $n$.  Also, if $|V_{n-1}| = t_{n-1}$ and
$|V_n|=t_n$ then  $E_n$ determines a $t_n \times t_{n-1}$ {\sl incidence
matrix}.  (See Figure 1 for an example.)

\bigskip

\psfig{figure=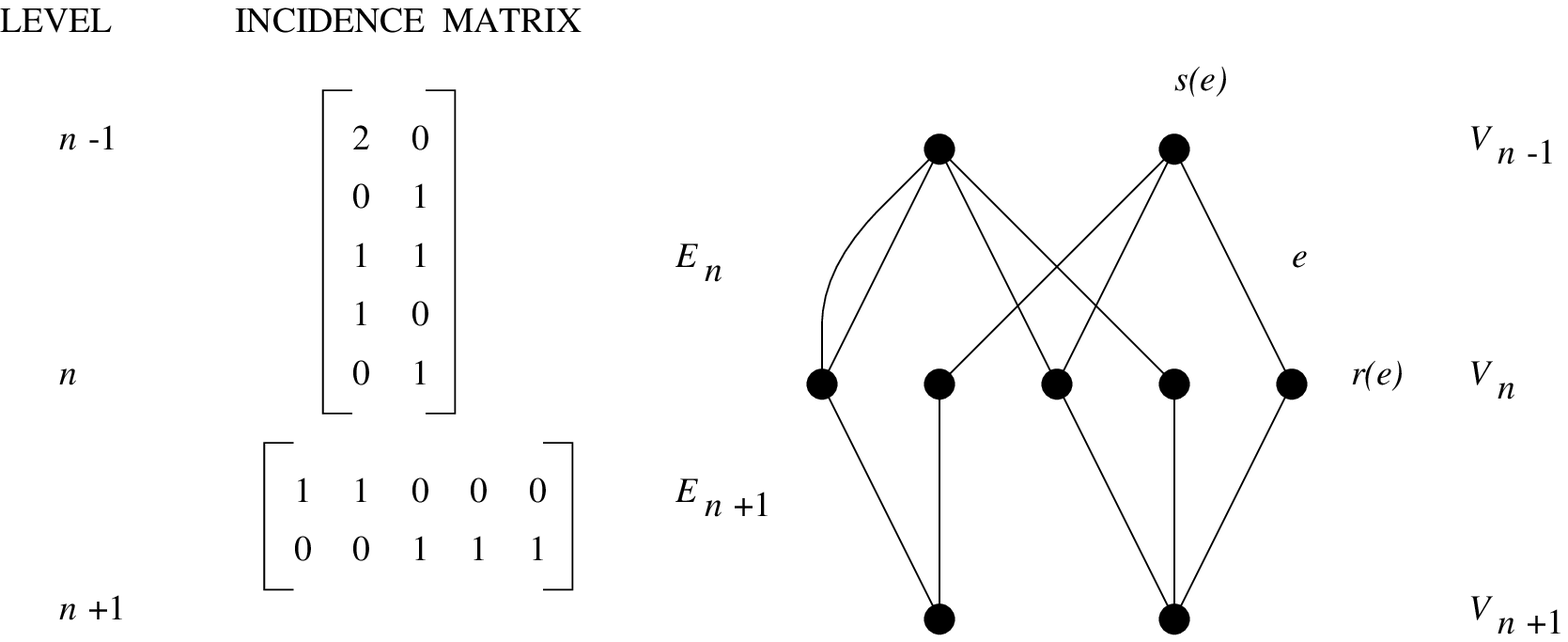,width=15cm,height=6cm}
\centerline{Figure 1.}

\bigskip

{\bf 1.1.2. Telescoping.}  Let $k , l \in {\Z}^+$ with $k < l$ 
and let $E_{k+1} \circ E_{k+2} \circ \cdots \circ E_l$ denote all paths 
from $V_k$ to $V_l$.
Specifically, 
$$ 
E_{k+1} \circ \cdots \circ E_l
$$
$$ = \{ ( e_{k+1} , \cdots , e_l ) |
e_i \in E_i , i = k + 1 , \cdots , l \ ; \
r(e_i ) = s(e_{i+1} ), i = k + 1 , \cdots , l - 1 \} . 
$$
We define $r((e_{k+1} , \cdots , e_l )) = r(e_l)$ and $s((e_{k+1} , 
\cdots , e_l )) \ = s(e_{k+1})  .$

Given a Bratteli diagram $(V,E)$ and a sequence 
$$
m_0 = 0 < m_1 < m_2 < \cdots 
$$
in ${\Z}^+$, we define the {\sl telescoping} of $(V,E)$ to $\{ m_n ; n\in \N \}$ 
as the new Bratteli diagram $(V' , E' )$, where  
$V_n ' = V_{m_{n}}$ and $E_n ' = E_{m_{n-1} + 1} \circ \cdots \circ E_{m_{n}}$ 
and the range and source maps are as above.

For example, if we remove level $n$ of Figure 1 we get a telescoping to 
levels $n-1$ and $n+1$ as indicated in Figure 2.  Note that the new
incidence matrix is the product of the two incidence matrices
of Figure 1.

\bigskip

\centerline{\psfig{figure=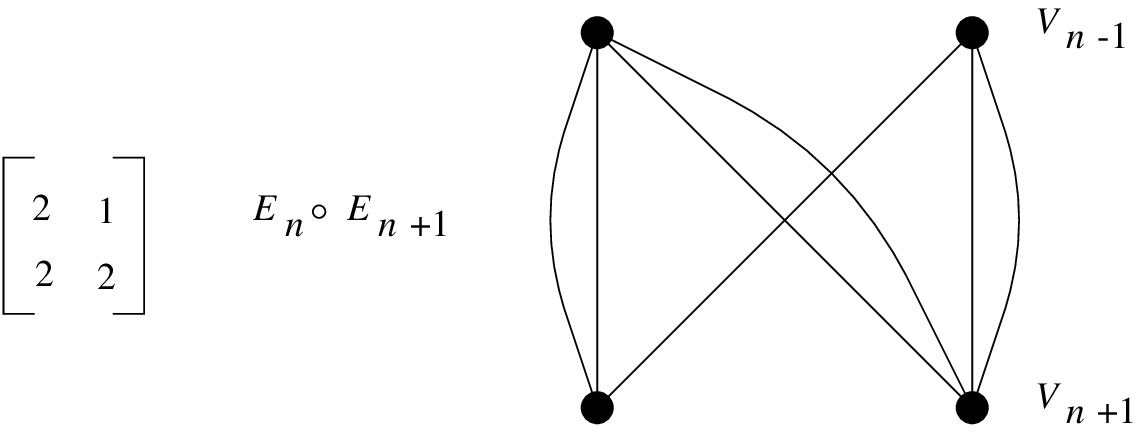,width=12cm,height=4cm}}
\centerline{Figure 2.}

\bigskip

We say that $(V,E)$ is a {\sl simple} Bratteli diagram if there exists a
telescoping $(V',E')$ of $(V,E)$ so that the incidence matrices of
$(V',E')$ have only non-zero entries at each level.

We let $\sim$ denote the equivalence relation on ordered Bratteli diagrams
generated by isomorphism and telescoping.  It is not hard to show that 
$(V^1,E^1) \sim (V^2,E^2)$ if and only if there exists a Bratteli diagram 
$(V,E)$ so that telescoping $(V,E)$ to odd levels $0 < 1 < 3 < 
\cdots$ yields a telescoping of either $(V^1,E^1)$ or $(V^2,E^2)$, and
telescoping 
$(V,E)$ to even levels $0 < 2 < 4 < \cdots$ yields a telescoping of the
other.  

\medskip

{\bf 1.2.  Ordered Bratteli diagrams and Bratteli-Vershik systems.}

\medskip

{\bf  1.2.1. Definition 2.}  {\sl An {\rm ordered Bratteli diagram} $(V,E, \geq )$
is a Bratteli
diagram $(V,E)$ together with a partial order $\geq$ on $E$ so that edges e,
e$'$
in $E$ are comparable if and only if $r(e) = r(e')$ ; in other words, we have
a linear order on each set $r^{-1} ( \{ v \} )$, where $v$ belongs to $V \setminus V_0$.}

\medskip

Note that if $(V, E, \geq )$ is an ordered Bratteli diagram and $k < l$ in
${\Z}^+$, then the set $E_{k+1} \circ E_{k+2} \circ \cdots \circ E_l$ of
paths from $V_k$ to $V_l$ may be given an induced (lexicographic) order as follows:
$$
(e_{k+1} , e_{k+2} , \cdots , e_l ) > (f_{k+1} , f_{k+2} , \cdots , f_l ) 
$$
if and only if for some $i$ with $k + 1 \leq i \leq l, \ e_j = f_j$
for $i < j \leq l$ and $e_i > f_i$.  It is a simple observation that if 
$(V, E, \geq )$ is an ordered Bratteli diagram and $( V' , E' )$ is a  
telescoping of $(V, E)$ as defined above, then with the induced order 
$\geq' , \; (V' , E' , \geq' )$ is again an ordered Bratteli diagram.  We
say that 
$(V' , E' , \geq' )$ is a {\sl telescoping} of $(V, E, \geq )$.

Again there is an obvious notion of isomorphism between ordered Bratteli
dia\-grams. 
We let $\approx$ denote the equivalence relation on ordered Bratteli diagrams
generated by isomorphism and by telescoping.  One can show that $B^1
\approx B^2$, where $B^1 = (V^1, E^1, \geq^1), \; B^2 = (V^2, E^2,
\geq^2)$, if and only if there exists an ordered Bratteli diagram $B
= (V,E, \geq)$ so that telescoping $B$ to odd levels $0 < 1 < 3 < 
\cdots$ yields a telescoping of either $B^1$ or $B^2$, and telescoping $B$
to even levels $0 < 2 < 4 < \cdots$ yields a telescoping of the
other.  This is analogous to the situation for the equivalence
relation $\sim$ on Bratteli diagrams as we discussed above.

\medskip

{\bf 1.2.2. Dynamics for ordered Bratteli diagrams.} Let $B = (V,E, \geq )$ be an ordered Bratteli diagram.  
Let $X_B$ denote the associated infinite path space, i.e.

$$ 
X_B = \{ (e_1 , e_2 , \cdots ) | e_i \in E_i , r (e_i ) = 
s (e_{i+1} )  ;  i = 1, 2, \cdots \}. 
$$

We will exclude trivial cases and assume henceforth that $X_B$ is an infinite
set.  Two paths in $X_B$ are said to be {\sl cofinal} if they have
the same tails, i.e. the edges agree from a certain stage. We topologize 
$X_B$ by postulating a basis of open sets, namely the family of cylinder sets

$$ 
U (e_1 , e_2 , \cdots , e_k ) = \{ ( f_1 , f_2 , \cdots ) \in
X_B | f_i = e_i , 1 \leq i \leq k \}. 
$$

Each $U (e_1 , \cdots , e_k )$ is also closed, as is  easily seen, and so we 
observe that $X_B$ becomes a compact Hausdorff space with a countable basis 
of clopen sets, i.e. a zero-dimensional space.  We call $X_B$ with this 
topology the {\sl Bratteli compactum} associated to $B = (V , E , \geq )$.
If $(V,E)$ is a simple Bratteli diagram, then $X_B$ has no isolated
points, and so is a {\sl Cantor set}.

\medskip

{\bf Notation.}  Let $x = (e_1, e_2, \cdots)$ be an element of $X_B$.  We
will call $e_n$ the $n$'th label of $x$ and denote it by $x(n)$.  We
let $X_B^{\max}$ (resp. $X_B^{\min}$) denote those elements $x$ of $X_B$
so that $x(n)$ is a maximal edge (resp. minimal edge) for each
$n=1,2, \cdots$. A simple argument shows that $X_B^{\max}$ (resp. 
$X_B^{\min}$) is non-empty.

\medskip

{\bf Definition 3.}  {\sl The ordered Bratteli diagram $B = (V,E , \geq)$ is
{\rm properly ordered} (called {\rm simple ordered} in
{\rm [HPS]}) if

\medskip

$\imath )$ $(V,E)$ is a simple Bratteli diagram.

$\imath\imath )$  $X_B^{\max}$, resp. $X_B^{\min}$, consists of only one
point $x_{\max}$, resp. $x_{\min}$.

}

\medskip

We can now define a minimal homeomorphism $V_B:X_B \rightarrow X_B$,
called the {\sl Vershik map} (or the {\sl lexicographic map}),
associated to the properly ordered Bratteli diagram $B=(V,E, \leq)$. 
We will call the resulting Cantor minimal system $(X_B, V_B)$ a
{\sl Bratteli-Vershik system}.

We let $V_B (x_{\max}) = x_{\min}$.  If $ x = (e_1, e_2, \cdots) \neq
x_{\max}$, let $k$ be the smallest number so that $e_k$ is not a
maximal edge.  Let $f_k$ be the successor of $e_k$ (and so $r(e_k) =
r(f_k))$.  Define $V_B (x) = y = (f_1, \cdots , f_{k-1}, f_k, e_{k+1}, 
e_{k+2} , \cdots)$, where $(f_1, \cdots, f_{k-1})$ is the
minimal edge in  $E_1 \circ E_2 \circ \cdots \circ E_{k-1}$ with 
range equal to $s(f_k)$.

It is a theorem that every Cantor minimal system is isomorphic to a
Bratteli-Vershik system.  A sketch of the proof of this result ("The
Bratteli-Vershik model theorem") will be presented in the next
subsection.  

\medskip

{\bf 1.3.  The Bratteli-Vershik model theorem -- Kakutani-Rohlin
partitions.}

\medskip

We shall sketch a proof of the following theorem using
Kakutani-Rohlin partitions.  We shall apply this later in the context
of substitution minimal systems.

\medskip

{\bf Theorem 5.} [HPS; Theorem 4.7] {\sl Let $(X,T,x)$ be a (pointed) Cantor
minimal system.  Then there exists a properly ordered Bratteli
diagram $B = (V,E, \geq)$ so that $(X,T,x)$ is pointedly isomorphic to $(X_B,
V_B, x_{\min})$, where $x_{\min}$ is the unique minimal path of $X_B$.}

\medskip
 
{\bf Definition 4.}  {\sl A {\rm Kakutani-Rohlin} partition of the
Cantor minimal system $(X,T)$ is a clopen partition ${\cal P}$ of the
kind:
$$ 
{\cal P} = \{ T^j Z_k | k \in A \ \ {\rm and} \ \ 0 \leq  j
< h_k \} 
$$
where $A$ is a finite set and $h_k$ is a positive integer.  The
$k'$th {\sl tower} of ${\cal P}$ is $\{ T^j Z_k | 0 \leq j < h_k \}$, and
the {\rm base} of ${\cal P}$ is the set $Z =
\mathrel{\mathop{\cup}\limits_{k \in A}} Z_k$.}

\medskip

Note that in [HPS] the roles of the base $Z$ and the top floors of
the towers are reversed.

\medskip

{\bf Sketch of proof of Theorem 5:}  Let $({\cal P}_n ; n\in \N)$ be a sequence
of Kakutani-Rohlin partitions with

$$ 
{\cal P}_n = \{ T^j Z_{n,k} | k \in A_n \ \ {\rm and} \ \
0 \leq j < h_{n,k} \}, \,\, {\cal P}_0 = \{ X \} 
$$

and with base $Z_n =  \mathrel{\mathop{\cup}\limits_{k \in A_n}}
Z_{n,k}$.  We say that this sequence is {\sl nested} if, for
each $n$,

\medskip

$\imath )$ $Z_{n+1} \subset Z_n.$

$\imath\imath )$ ${\cal P}_{n+1} \succ {\cal P}_n$ as partitions.

\medskip

To the nested sequence $({\cal P}_n ; n\in \N )$ we associate an ordered
Bratteli diagram $B = (V, E, \geq)$ as follows:

The $|A_n|$ towers in ${\cal P}_n$ are in 1-1 correspondence with
$V_n$, the set of vertices at level $n$.  Let $v_{n,k} \in
V_n$ correspond to the tower $S_{n,k} = \{ T^j Z_{n,k} | 0 \leq j <
h_{n,k} \}$ in ${\cal P}_n$.  $S_{n,k}$ will traverse towers in
${\cal P}_{n-1}$ in a certain order, say $S_{n-1, i_1}, \cdots ,
S_{n-1, i_m}$.  (Thus $S_{n,k}$ is partitioned into $m$ parts.)  We
associate $m$ edges, ordered as $e_{1,k} < e_{2,k} < \cdots < 
e_{m,k}$, with $r(e_{j,k} ) = v_{n,k}$ and $s(e_{j,k}) = v_{n-1, i_j}$.  
$E_n$ is the disjoint union over $k \in A_n$ of the edges with range
equal to $v_{n,k}$.

Suppose moreover that:

\medskip

$\imath\imath\imath )$ The intersection of the bases $(Z_n ; n\in \N)$ of the partitions
$({\cal P}_n ; n\in \N)$ consists of one point only, say $x$.

$\imath \nu )$ The sequence of partitions spans the topology of $X$.

\medskip

Then the systems $(X,T,x)$ and $(X_B, V_B, x_{\min})$ are pointedly
isomorphic, where $B = (V,E, \geq)$ is a properly ordered Bratteli
diagram with unique minimal path $x_{\min}$.  In fact, the isomorphism $F: X \rightarrow X_B$ is defined as follows:  For $y \in
X$, the path $F(y)$ passes through the vertex in $V_n$ that
corresponds to the tower in ${\cal P}_n$ where $y$ is located.  Say
$F(y)$ passes through the vertices $w$ of $V_{n-1}$ and $v$ of $V_n$, corresponding to the tower $S_w$ in ${\cal
P}_{n-1},\;S_v$ in ${\cal P}_n$, respectively.  Then $F(y)(n)$ is the
$i'$th edge among the ordered edges $e$ with $r(e) = v, \; y$ being
"picked up" the $i'$th time $S_v$ traverses one of the towers at level
$n-1$.  (Necessarily, $S_v$ will traverse the tower $S_w$ the $i'$th
time.)  It is easily seen that $F(x) = x_{\min}$ by this correspondence.

Finally, one obtains Kakutani-Rohlin partitions $({\cal P}_n ; n\in \N)$ with
the properties listed above by choosing a nested sequence $(Z_n ; n\in \N)$ of clopen
sets shrinking to $x$, i.e. $Z_n \supseteq Z_{n+1}$ and 
$\mathrel{\mathop{\cap}\limits_{n}} Z_n = \{ x \}$.  The partition
${\cal P}_n$ is obtained by building towers over $Z_n$ by considering
the return map to $Z_n$ (For details, cf. [HPS]). 
\cqfd

\medskip

{\bf 1.4. Dimension groups.}

\medskip

{\bf 1.4.1. Definition 5.}  A {\sl dimension group} is an ordered, countable,
torsion free abelian group $G$ which is unperforated and satisfies the Riesz
interpolation property.  Specifically, let $G^+$ denote the positive cone of
$G$ and let $a \leq b$ denote $b - a \in G^+$.  Then

\medskip

$\imath )$ $G^+ + G^+ \subseteq G^+$

$\imath\imath )$ $G^+ - G^+ = G$

$\imath\imath\imath )$ $G^+ \cap (-G^+) = \{ 0 \}$

$\imath \nu )$ If $a \in G$ and $na \in G^+$ for some $n \in {\bf N}$, then $a
\in G^+$ (unperforation).

$\nu )$ If $a_1, a_2, b_1, b_2 \in G$ with $a_i \leq b_j \;(i,j = 1,2)$,
there exists $c \in G$ with $a_i \leq c \leq b_j$ (Riesz interpolation
property).

\medskip

We say that $(G, G^+)$ is {\sl simple} if $G$ has no non-trivial order
ideals $J$, i.e. $J$ is a proper subgroup so that $J = J^+ - J^+$ (where
$J^+ = J \cap G^+)$ and so that if $b$ is an element of $J$ and $a$ an element of $G$ with $0 \leq a \leq b $, then $a$ belongs to $J$.

We say that an element $u$ of $G^+ \setminus \{ 0 \}$ is an {\sl order unit} if $G^+ =
\{ a \in G| 0 \leq a \leq nu$ for some $n \in {\bf N}\}$.  (Observe that for
$G$ simple any $u \in G^+ \setminus \{ 0 \}$ is an order unit.)

To the Bratteli diagram $(V,E)$ is associated a dimension group which we
denote by $K_0 (V,E)$ --- the notation is motivated by the connection to
$K$-theory, see below.  In fact, to the Bratteli diagram $(V,E)$ is
associated a system of ordered groups and order-preserving homomorphisms 

$$
\Z^{|V_0|}\fleche{\varphi_1}
\Z^{|V_1|}\fleche{\varphi_2}
\Z^{|V_2|}\fleche{\varphi_3}
\Z^{|V_3|}\fleche{\varphi_4}\cdots
$$

where $\varphi_n$ is given by matrix multiplication with the incidence
matrix between levels $n-1$ and $n$ of the Bratteli diagram.  By definition
$K_0 (V,E)$ is the inductive limit of the system above endowed with the
induced order.  $K_0 (V,E)$ has a distinguished order unit, namely the
element of $K_0 (V,E)^+$ corresponding to the element $1 \in {\Z}^{|V_0|}
= {\Z}$.  One can show that $(V,E) \sim (V',E')$ if and only if
$K_0(V,E)$ is order isomorphic to $K_0 (V',E')$ by a map sending the
distinguished order unit of $K_0 (V,E)$ to the distinguished order unit of
$K_0(V',E')$. The dimension group $K_0 (V,E)$ associated to $(V,E)$ is
simple if and only if $(V,E)$ is simple.

\medskip

{\bf Remark.}  It is a theorem that all dimension groups arise from Bratteli
diagrams as described above [EHS].

\medskip

{\bf 1.4.2.}  We now introduce the definition which will relate Cantor
minimal systems to dimension groups.

\medskip

{\bf Definition 6.}  {\sl Let $(X,T)$ be a Cantor minimal system.  Let $C(X, {\bf
Z})$ denote the continuous functions on $X$ with values in ${\Z}$ --- so
$C(X, {\Z})$ is a countable abelian group under addition.  Let 
$$ 
K^0 (X,T) = C (X, {\Z})\slash \partial_T C(X, {\Z}) 
$$
where $\partial_T: C(X,{\Z}) \rightarrow C(X,{\Z})$ denotes the {\rm
coboundary operator} $\partial_T (f) = f - f \circ T$, and $f - f \circ T$ is
called a {\rm coboundary}.  Define the {\rm positive cone} 
$$ 
K^0 (X,T)^+ = \{ [f] | f \in C (X, {\Z}^+) \} 
$$
where $[\, \cdot \,]$ denotes the quotient map and ${\Z}^+ = \{ 0,1,2,
\cdots \}$.  $K^0 (X,T)$ has a {\rm distinguished order unit}, namely [1] = ${\bf
1}$, where 1 denotes the constant function one.
}

\medskip

{\bf Theorem 6.} [HPS; Theorem 5.4 and Corollary 6.3]  {\sl Let $(X,T)$ be a Cantor minimal
system.  Let $B = (V,E, \geq)$ be the associated properly ordered Bratteli
diagram (having chosen a base point in $X$, cf. Theorem  5).  Then
$$ 
K^0 (X,T) = K_0 (V,E) 
$$

as ordered groups with distinguished order units.  Furthermore, every
simple dimension group $G \; (G \neq {\Z})$ arises in this manner.}

\medskip

{\bf Remarks.}  

$\imath )$ 
In [GPS] it is shown that $K^0 (X,T)$, as an ordered group with
distinguished order unit, is a complete invariant for strong orbit
equivalence of Cantor minimal systems.

$\imath\imath )$ 
$K^0 (X,T)$ is order isomorphic, by a map preserving the
distinguished order units, to the $K_0$-group of the $C^\ast$-crossed
product associated to $(X,T)$.

$\imath\imath\imath )$
The group $K^0 (X,T)$, as an abstract group without order, is
isomorphic to the first \u{C}ech cohomology group $H^1 (\hat{X}, {\Z})$
of the suspension $\hat{X}$ of $(X,T)$, where $\hat{X}$ is obtained from $X
\times [0,1]$ by identifying $(x,1)$ and $(Tx, 0)$.

\medskip

{\bf 1.4.3.}  We will describe an alternative method to associate a
dimension group to a nested sequence of Kakutani-Rohlin partitions (cf. subsection 1.3)
without invoking Bratteli diagrams explicitly.  This method, close to the
approach taken in [GW], will later be applied to substitution miminal
systems.  (Incidentally, the proof of the first part of Theorem 6
is an immediate consequence of the lemma below.)
So let $({\cal P}_n ; n\in \N)$ be a nested sequence of Kakutani-Rohlin partitions
(not necessarily satisfying conditions $\imath\imath\imath)$ and $\imath\nu)$ of 1.3), with

$$ 
{\cal P}_n = \{ T^j {Z}_{n,k} | k \in A_n , 0
\leq j < h_{n,k} \}, \;\;  {\cal P}_0 = \{  X \} 
$$

and with base $Z_n = \cup_{k \in A_n} Z_{n,k}$.  For each $n$, let $C_n$ be
the subgroup of $C(X, {\Z})$ consisting of functions which are constant
on each element of the partition ${\cal P}_n$, and let $C^+_n = C_n \cap  C
(X, {\Z}^+)$.  Let $H_n$ be the subgroup of $C_n$ consisting of the
functions $f \in C_n$ which have a null sum over each tower of ${\cal P}_n$,
i.e. such that

$$ 
\sum^{h_{n,k}-1}_{j=0} f(T^j x) = 0 \ \ \forall  k \in A_n ,
\ \ \forall x \in  Z_{n,k}. 
$$

Let $K_n$ be the quotient group $C_n/H_n, \, K^+_n$ the projection of
$C^+_n$ in this quotient, and ${\bf 1}_n$ the projection of the constant
function $1$.  $(K_n, K^+_n, {\bf 1}_n)$ is an ordered group with order
unit.  Obviously, $C_n \subset C_{n+1}, \; C^+_n \subset C^+_{n+1}$ and $H_n
\subset H_{n+1}$; these inclusions induce a morphism $j_{n+1}: K_n
\longrightarrow K_{n+1}$ of ordered groups with order units.  The direct
limit $(K({\cal P}), K({\cal P})^+, {\bf 1})$ of the sequence
$$
K_0\fleche{j_1}
K_1\fleche{j_2}
K_2\fleche{j_3}
\cdots
\fleche{j_n}
K_n\fleche{j_{n+1}}
K_{n+1}\fleche{j_{n+2}}
\cdots
$$
of ordered groups with order units is called {\sl the dimension group of the
sequence of partitions}, and written $K({\cal P}_n; n \in \N)$. 
It is easy to check that every element of $H_n$ is a coboundary.  Thus the
inclusions $C_n \hookrightarrow C (X, {\Z})$ induce a morphism $\lambda:
K ({\cal P}_n ; n \in \N) \longrightarrow K^0 (X,T)$ of ordered groups with
order units.

\medskip

{\bf Lemma 7.}  {\sl As an ordered group with order unit, $K({\cal P}_n ; n \in \N )$ is isomorphic to the dimension group $K_0 (V,E)$, where $B = (V,E,
\geq)$ is associated to the sequence $({\cal P}_n ; n\in \N)$ of partitions as
described in subsection 1.3 (disregard the ordering of the edges).}

\medskip

{\bf Proof.}  For each $f$ of $C_n$, let $\gamma_n (f)$ be the vector of ${\Z}^{|V_n|}$ obtained by summing the values of $f$ over each tower of
$P_n$:  For each $v_{n,k}$ belonging to $V_n$ (corresponding to the tower with base
$Z_{n,k})$, let

$$ 
(\gamma_n (f))_k = \sum^{h_{n,k}-1}_{j=0} f(T^j x) 
$$

for an arbitrary point $x \in Z_{n,k}$.  Now $\gamma_n$ is a positive
homomorphism from $C_n$ onto ${\Z}^{|V_n|}$, it maps $C^+_n$ onto $({\Z}^{|V_n|})^+$, and its kernel is $H_n$.  Thus we can
identify the ordered groups $K_n = C_n/H_n$ and ${\Z}^{|V_n|}$ for
$n=0,1,2,\cdots . \; $  By these identifications the homomorphism $j_{n+1} :
K_n \longrightarrow K_{n+1}$ corresponds to the homomorphism ${\bf
Z}^{|V_n|} \longrightarrow {\Z}^{|V_{n+1}|}$ that we get from the
Bratteli diagram $(V,E)$.  Therefore the two inductive limits are order
isomorphic.  Furthermore, it is easy to verify that the distinguished order
units are mapped to each other by this isomorphism. 
\cqfd

\medskip

{\bf 1.5. Stationary diagrams and dimension groups.}

\medskip

{\bf 1.5.1. Definition 7.}  A Bratteli diagram $(V,E)$ is {\sl stationary} if
$k= |V_1|=|V_2|= \cdots$ and if (by an appropriate labeling of the vertices)
the incidence matrices between level $n$ and $n+1$ are the same $k \times k$
matrix $C$ for all $n=1,2,\cdots. \;$  In other words, beyond level $1$ the
diagram repeats.  (Clearly we may label the vertices in $V_n$ as $V(n,a_1),
\cdots , V(n, a_k)$, where $A=\{ a_1, \cdots, a_k\}$ is a set of $k$
distinct symbols.)

$B = (V,E, \geq)$ is a {\sl stationary ordered Bratteli diagram} if $(V,E)$
is stationary, and the ordering on the edges with range $V(n,a_i)$ is the
same as the ordering on the edges with range $V(m,a_i)$ for $m,n=2,3,\cdots$
and $i=1,\cdots, k$.  In other words, beyond level $1$ the diagram with the
ordering repeats.  (For each $a_i$ in $A=\{a_1, \cdots , a_k\}$ and each
$n=2,3,\cdots,$ we thus get an ordered list of edges whose range is
$V(n,a_i)$.  By the stationarity of the ordering of $B$ we thus get a well
defined map from $A$ to $A^{+}$ (the set of non-empty words on $A$), 
by taking the sources of the edges in question.)

$(G,G^+)$ is a {\sl stationary dimension group} if $G$ is order isomorphic
to $K_0(V,E)$, where $(V,E)$ is a stationary Bratteli diagram.  $K_0(V,E)$
is completely determined by the incidence matrix $C$ of $(V,E)$ --- we
disregard the distinguished order unit.  Also, $K_0 (V,E)$ is simple if and only if $C$
is a {\sl primitive} matrix, i.e. a certain power of $C$ has only non-zero
entries.

\medskip

In this paper we will only encounter stationary Bratteli diagrams.  In
Figure 3 we exhibit two examples of stationary ordered Bratteli diagrams;
the one on the left is properly ordered, while the one on the right is not properly ordered (having in fact two max and two min paths).
The dimension group associated to the two diagrams (strip the order
structure) is ${\Z}[1/2]$, the dyadic rationals, with obvious ordering and
with distinguished order unit equal to $1$.

\bigskip

\centerline{\psfig{figure=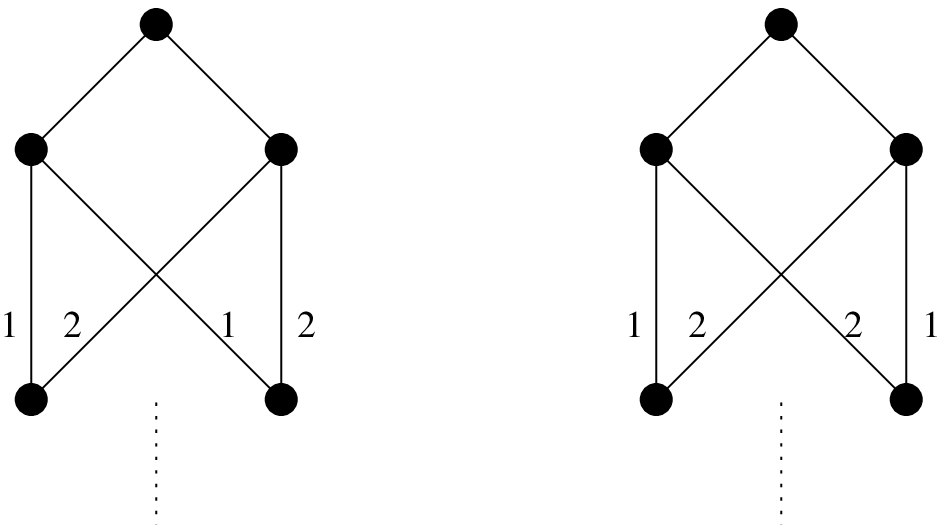,width=10cm,height=5cm}}
\centerline{Figure 3.}
\bigskip

{\bf 1.5.2.}  The question of the unicity of the Bratteli-Vershik model is
addressed in the following proposition.

\medskip

{\bf Proposition 8.} [HPS; Theorem 4.7]  {\sl Let $B^i = (V^i,E^i, \geq^i)$ be
associated to the pointed Cantor minimal system $(X_i, T_i, x_i), i=1,2$,
according to Theorem 5.  Then $(X_1, T_1, x_1)$ is pointedly
isomorphic to $(X_2, T_2, x_2)$ if and only if $B^1 \approx B^2$ (cf. part 1.2.1).}

\medskip

By a modification of the "symbol splitting" procedure in [GPS; Sect. 3],
Forrest proved the following lemma, which will be useful to us.

\medskip

{\bf Lemma 9.}  [Fo; Lemma 15]  {\sl Let $B=(V,E,\geq)$ be a stationary, properly
ordered Bratteli diagram.  Then $B \approx B' = (V',E', \geq')$, where $B'$
is again stationary and properly ordered, and with the added property that
there are no multiple edges between level $0$ and level $1$.  (By the above
proposition, $(X_B, V_B)$ and $(X_{B'}, V_{B'})$ are isomorphic.)}

\medskip

{\bf Proof.}  By taking a sufficiently high power of the incidence matrix
$C$ of $(V,E)$ (this corresponds to a periodic telescoping of the diagram),
we may assume that each row sum of $C$ is greater or equal to the maximum
number of edges between a vertex at level 1 and the top vertex (i.e. level
0).  We now proceed by introducing new vertices between two successive
levels, the number of which is the same as $|E_1|$, i.e. the number of edges
between level 0 and level 1.  One may now construct a properly ordered
Bratteli diagram with these added levels, so that one gets the original by
telescoping.  (This construction is not unique.)  By instead telescoping to
the new levels introduced one gets the desired $B' = (V',E', \geq')$.
We will give a diagrammatic example which will illustrate how to proceed in
the general case (cf. Figure 4). 
\cqfd

\bigskip

\centerline{\psfig{figure=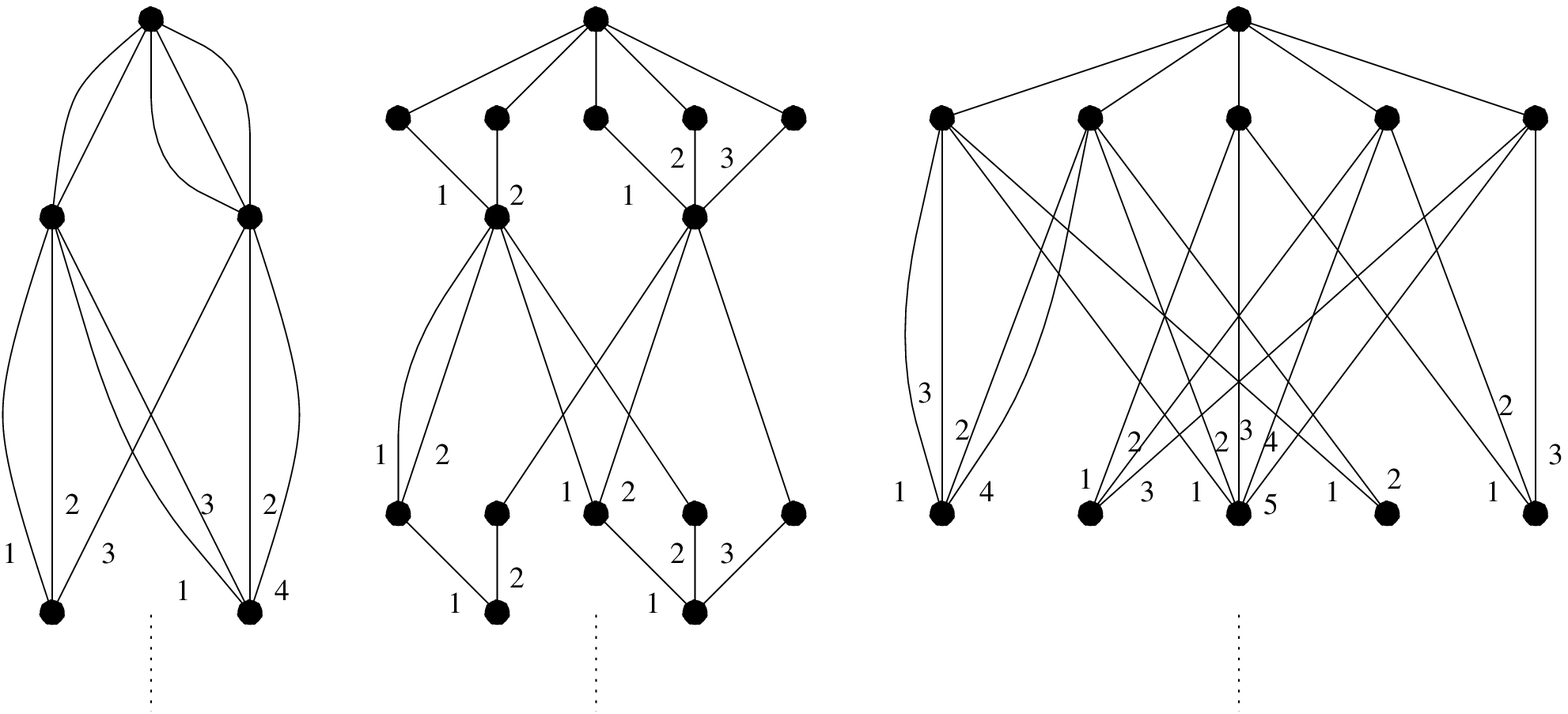,width=16cm,height=7cm}}

\centerline{Figure 4.}

\bigskip

{\bf 1.6. Kakutani equivalence.}  

\medskip

{\bf Definition 8.} {\sl The Cantor minimal systems $(X,T)$ and $(Y,S)$ are {\rm
Kakutani equivalent} if they have (up to isomorphism) a common derivative,
i.e. there exist clopen sets $U$ (in $X$) and $V$ (in $Y$), respectively, so
that the induced systems on $U$ and $V$, respectively, are isomorphic.}

\medskip

We will relate Kakutani equivalence to Bratteli diagrams --- the relevant
fact being change of the order unit.

Observe first that if $(V,E)$ is a Bratteli diagram with associated
dimension group $G=K_0 (V,E)$, then any finite change of $(V,E)$, i.e.
adding and/or removing a finite number of edges (vertices), thus changing
$(V,E)$ into a new Bratteli diagram$(V',E')$, does not change the
isomorphism class of $G$, but does change the order unit.  In fact,
$G'=K_0(V',E')$ is order isomorphic to $G$, but the distinguished order
units are not necessarily preserved by the isomorphism.  Clearly, any change
of order unit of $G$ may be obtained by such a procedure.

Likewise, if $B=(V,E,\geq)$ is a properly ordered Bratteli diagram we may
change $B$ into a new properly ordered Bratteli diagram $B'=(V',E',\geq')$
by making a finite change, i.e. adding and/or removing any finite number of
edges (vertices), and then making arbitrary choices of linear orderings of
the edges meeting at the same vertex (for a finite number of vertices).  So
$B$ and $B'$ are cofinally identical, i.e. they only differ on finite
initial portions.  (Observe that this defines an equivalence relation on the
family of properly ordered Bratteli diagrams.)

\medskip

{\bf Theorem 10.} [GPS; Theorem 3.8] {\sl Let $(X_B,V_B)$ be the Bratteli-Vershik
system associated to the properly ordered Bratteli diagram $B=(V,E,\geq)$.  
Then the Cantor system $(Z,\psi)$ is Kakutani equivalent to $(X_B,V_B)$ if
and only if $( Z,\psi)$ is isomorphic to $(X_{B'},V_{B'})$, where $B' = (V',E',
\geq')$ is obtained from $B$ by a finite change as described above.}

\medskip

We indicate briefly the ingredients of the proof:  Every clopen set of $X_B$
is a finite union of cylinder sets.  By telescoping appropriately one may
assume that the clopen sets in question are disjoint unions of cylinder sets
of the form $U(e)$, where $e \in E_1$ (cf. 1.2.2.).  This is achieved by
making a finite change of the original diagram.  One may thus assume at the
outset that the only changes that are done to the diagram occur between
level $0$ and level $1$.  With this simplification the proof follows by
observing that removing edges corresponds to inducing on a clopen set.

%% file: section2.tex
{\gros 2. Basic facts about substitution dynamical systems.}

\bigskip

{\bf 2.1. Notations: words, sequences, morphisms.} 

\medskip

 An {\sl alphabet} is a finite set of symbols called
 {\sl letters}. If $A$ is an alphabet, a {\sl word } on $A$ is a
finite (non empty) sequence of letters; $A^+$ is the set of words.
For $u=u_1u_2\ldots u_n\in A^+$, $|u|=n$ is the {\sl length} of $u$; for each letter $a$, 
we write $|u|_a$ the number of occurrences of $a$
in $u$; the vector $(|u|_a\semi a\in A)$ is sometimes called the {\sl composition vector}
of $u$.
$A^*$  consists of $A^+$ and the empty word $\emptyset$ of length $0$.
. 

Given a word $u=u_1\ldots u_m$ and an interval $J=\{i,\ldots,j\}$ contained in
$\{1,\ldots,m\}$, we write $u_J$ to denote the word  $u_iu_{i+1}\ldots u_j$. We extend this notation in an obvious way to infinite intervals. A {\sl factor} of $v$
is a word $u$ such that $u=v_J$ for some interval $J\subset\{1,\ldots,m\}$; we write
$u\prec v$.

Elements of $A^\Z$ are called {\sl sequences} over the alphabet $A$. For a sequence $x$
we use the notation $x_J$ and the term {\sl factor} exactly as for a word; the {\sl
language} $\L(x)$ of the sequence $x$ is the set of words which are factors of $x$.

Let $A,B$ be two alphabets, and $\phi\colon A\to B^+$ a map. By
concatenation, $\phi$ can be extended to a map $A^+\to B^+$ and
also to a map $A^\Z\to B^\Z$. As a general rule, these maps will be
denoted by $\phi$ also. 

The following definitions will be used in Section 4:

\medskip
{\bf Definition 9. } {\sl Let $A$ be an alphabet, and $\RR$ a finite subset of
$A^+$.

$\imath$) $\RR$ is a code if every word $u\in A^+$ admits at most one
decomposition in a concatenation of elements of $\RR$.

$\imath\imath$) $\RR$ is a circular code if it is a code and moreover:

if
 $$
w_1,\ldots,w_j,\   w, \ w'_1 , \ldots , w'_k \in \RR \ ; \ 
s\in A^+ \hbox{ and } t\in A^*
$$
are such that
 $$w=ts\hbox{ and }w_1\ldots w_j=sw'_1\ldots w'_kt$$
then $t$ is the empty word.}

(it follows that $j=k+1$, $w_{i+1}=w'_{i}$ for $1\leq i\leq k$ and
that $w_1=s$.)
\medskip

Circular codes have a property of unique decomposition of
sequences: Suppose that $\RR$ is a circular code on the alphabet $A$, and that some
$x\in A^\Z$ can be decomposed in a the concatenation of the words
$(w_k\semi k\in\Z)$ belonging to $\RR$, i.e. that 
$$
x=\ldots w_{-3}w_{-2}w_{-1} \mid w_0w_1w_2\ldots
$$
where the vertical bar separates $x_{(-\infty,-1]}$ from $x_{[0,+\infty)}$.
Then this decomposition is unique. Moreover, for $n\in \Z$, the shifted sequence $T^nx$ 
(given by $(T^nx)_k=x_{n+k}$ for all $k$) can
be decomposed in a concatenation of elements of $\RR$ if and only if $n$ is
one of the numbers:
$$
-|w_rw_{r+1}\ldots w_{-1}| \ (r<0)\semi 0\semi
|w_0w_1\ldots w_r|\ (r \geq 0)\ .
$$

\medskip

{\bf 2.2. Subshifts.}

\medskip

For every alphabet $A$ we denote by $T$ the shift on $A^\Z$, except when some
ambiguity can occur; $T$ is defined by
$$(Tx)_n=x_{n+1}\hbox{ for every }x\in A^\Z\hbox{ and every }n\in \Z\ .$$
$A$ is endowed with the discrete topology, and $A^\Z$ with the product topology; thus
$A^\Z$ is a compact metric space, and $T$ a homeomorphism of this space. 

When $B$ is an alphabet and $\phi: A\to B^+$ a map, the corresponding map
$\phi:A^\Z\to B^\Z$ is continuous, and satisfies
$$\phi(Tx)=T^{|\phi(x_0)|}\phi(x)\ .$$

A {\sl subshift} $(X,T)$ on the
alphabet $A$ is a closed $T$-invariant subset $X$ of $A^\Z$, endowed with the restriction
of $T$ to $X$, we denote it $T$ too. Given a subshift $(X,T)$ on the alphabet $A$, and two
words $u,v$ on $A$, we write: 
$$
[u]=\{x\in X\semi x_{[0,|u|)}=u\} \hbox{ and } [u.v]=\{x\in
X\semi x_{[-|u|,|v|)}=uv\}\ ;
$$
subsets of $X$ of this kind are called {\sl cylinder sets}.
Remark that these notations depend on the given subshift, which is supposed to be determined
by the context.

The {\sl language} $\L (X)$ of the subshift $(X,T)$ is the set of words which are factors of
at least one element of $X$,  and we have
$$X=\{x\in A^\Z\semi \L(x)\subset \L(X)\}\ .$$
The subshift $(X,T)$ is minimal if and only if $\L(x)=\L(X)$ for all $x\in X$.

The {\sl subshift spanned } by a sequence $x\in A^\Z$ is the closure $\barO_x$ of the
$T$-orbit of $x$ endowed with the restriction of the shift to $\barO_x$; it is
characterized by $\L(\barO_x)=\L(x)$.  This subshift is minimal if and only if $x$ is {\sl
uniformly recurrent} i.e.:
$$
(\forall u\in\L(x))(\exists n\geq 1)(\forall v\in\L(x))
(|v|\geq n\Rightarrow u\prec v ).
$$
\medskip

{\bf 2.3. Definition of substitution dynamical systems.}

\medskip

{\bf 2.3.1.} A {\sl substitution} on the alphabet $A$ is a map $\sigma\colon
A\to A^+$. Using the extension to words by concatenation, $\sigma$
can be iterated; for each integer $n>0$, $\sigma^n\colon A\to A^+$
is again a substitution.

In this paper we only consider {\sl primitive} substitutions, i.e. 
substitutions $\sigma$ on $A$ such that:
$$
(\exists n>0)(\forall a,b\in A)(b\prec \sigma^n(a))\;\; {\rm and} \;\; (\exists a\in A)(\lim_{n\rightarrow +\infty} |\sigma^n (a)| = +\infty).
$$

We denote by $\L(\sigma)$  {\sl the  language of} $\sigma$, i.e. the set of
words on $A$ which are factors of $\sigma^n(a)$ for some $a\in A$ and
some $n\geq 1$, and $\xs$ the subshift of $A^\Z$ associated to this
language, i.e. the
 set of $x\in A^\Z$ whose every finite factor belongs to
$\L(\sigma)$; $\xs$ is closed in $A^\Z$, and invariant under the
shift; we denote by $\ts$  the restriction of the shift to $\xs$. The
dynamical system $(\xs,\ts)$ is called {\sl the substitution
dynamical system } associated to $\sigma$. It is classical [Qu] that

\medskip
{\sl Every substitution dynamical system is minimal and
uniquely ergodic.}

\medskip
{\bf 2.3.2.} In the literature, substitution dynamical systems are often
defined by a different (but equivalent) method, using  fixed
points:

For every integer $p>0$, the substitution $\sigma^p$ defines the
same language, thus the same system, as $\sigma$ does.
Substituting $\sigma^p$ for
$\sigma$ if needed,  we can assume that there exist two letters
$r,\ell\in A$ such that:

{\parindent = 2cm 
$\imath$) $r$ is the last letter of $\sigma(r)$;

$\imath\imath$) $\ell$ is the first letter of $\sigma(\ell)$

$\imath\imath\imath$)  $r\ell\in\L(\sigma)$.
\par}

\medskip

Whenever $r$ and $\ell$ satisfy the conditions $\imath$) and 
$\imath\imath$), 
it is easy to
check that there exists a unique $\omega\in A^\Z$ such that
$$\omega_{-1}=r\semi \omega_0=\ell\hbox{ and }
\sigma(\omega)=\omega\ .$$
Such an $\omega$ is called {\sl a  fixed point of} $\sigma$. If $r$ and
$\ell$ satisfy also $\imath\imath\imath$), we say that $\omega$ is an
{\sl admissible fixed point of $\sigma$}.

If $\omega$ is an admissible fixed point of $\sigma$, then
 $\xs$ is the closure of the orbit of $\omega$ for the shift. This property
is often taken of the definition of $\xs$ in the literature. More precisely,
we have:

\medskip
{\sl 
Given a fixed point $\omega$ of the primitive substitution $\sigma$,  the
following conditions are equivalent:

$\imath$) $\omega$ is admissible.

$\imath\imath$) $\omega\in \xs$.

$\imath\imath\imath$) $\xs$ is the subshift spanned by $\omega$.

$\imath\nu$) $\omega$ is uniformly recurrent.
}

\medskip

In all this paper, we shall use the following convention: when we say that
$\omega$ is a fixed point of $\sigma$, we mean that $\omega$ is a fixed
point of $\sigma^p$ for some $p$; as $\sigma$ and $\sigma^p$ define the
same system, this convention can't lead to any misunderstanding; moreover,
with this convention, every substitution has at least one admissible fixed
point.
 
\medskip

{\bf 2.3.3.} Some technical difficulties in the study of substitution
dynamical system arise from the fact that a given substitution can have
several admissible fixed points. We now introduce a class of substitutions
which are easier to study:

\medskip

{\bf Definition 10.} {\sl A substitution $\sigma$ on the alphabet $A$ is
{\rm proper} if there exists an integer $p>0$ and two letters $r,\ell\in A$ such
that:

{\parindent 2cm
$\imath$) For every $a\in A$, $r$ is the last letter of $\sigma^p(a)$;

$\imath\imath$) For every $a\in A$, $\ell$ is the first letter of $\sigma^p(a)$
\par}
}
\medskip
A proper substitution has only one fixed point. We shall see later (in
Section 5) that every substitution dynamical system is isomorphic to the
system associated to some proper substitution, which can be explicitly
constructed.
\medskip

{\bf 2.4. Structure of substitution dynamical systems.}
\medskip

In the sequel, $(\xs,\ts)$ is the system associated to the 
primitive substitution $\sigma$  on the alphabet $A$, 
$\omega$ is an admissible fixed point of $\sigma$, $r=\omega_{-1}$ and
$\ell=\omega_0$.

There exist substitutions $\sigma$ such that the system $\xs$ is finite, i.e.
that every $x$ in $\xs$ is periodic, or equivalently, that
 $\omega$ is periodic. As these substitutions are of
little interest from a dynamical point of view, we consider henceforth only here
{\sl aperiodic} substitutions, i.e. substitutions giving rise to
infinite systems. Note that there is an algorithm [Pa,HL] which
decides whether a given substitution is aperiodic or not.

\medskip

The systems arising from (primitive, aperiodic) substitutions have
a simple self--similar structure we explain now. Remark first that  by the first
definition of $\xs$ we have $\sigma(\xs)\subset\xs$.
The result below  says that, although
neither $\sigma: A\to A^+$ is  assumed to be one to one, nor $\{\sigma(a)\semi a\in A\}$
assumed to be a code, the restriction of $\sigma$ to $\xs$ behaves almost as these
properties were true;  let's begin with some notation.

As $\omega=\sigma(\omega)$, it can be written
$$\omega=\ldots\sigma(\omega_{-2})\sigma(\omega_{-1})\mid 
\sigma(\omega_0)\sigma(\omega_1)\ldots$$
where the vertical bar separates $\omega_{(-\infty,-1]}$ and 
$\omega_{[0,+\infty)}$. Let $E=\{ e_j\semi j\in\Z\}$ be the set of ``natural cutting
points'' arising from this decomposition, defined by
$$e_j=\left\{
\matrix{
-\bigl|\sigma\bigl(\omega_{[j,0)}\bigr)\bigr|& \hbox{ if } & j<0\cr
0 & \hbox{ if } & j=0\cr
\bigl|\sigma\bigl(\omega_{[0,j)}\bigr)\bigr| &\hbox{ if } & j>0\cr
} \right.
$$
\medskip

{\bf Theorem 11.} [Mo1,Mo2] {\sl Let $\sigma	$ be an aperiodic primitive substitution on the alphabet
$A$.  

$\imath$) {\rm [Mo1]} There exists $L>0$ such that
$$n\in E,\; m\in\Z,\; \omega_{[n-L,n+L)}=\omega_{[m-L,m+L)}\ \Rightarrow m\in E \ .$$
$\imath\imath$) {\rm [Mo2]} There exists $M>0$ such that 
$$i,j\in \Z,\; \omega_{\textstyle [e_i-M,e_i+M)}=
\omega_{\textstyle [e_j-M,e_j+M)}\ \Rightarrow 
\omega_i=\omega_j\ .$$
}
\medskip

{\bf Remark.} Earlier proofs of the same (or similar) results do
exist [Ma], but are false. Until Mosse's papers, people working
with substitutions had to make extra hypothesis (`recognizability',
injectivity) on the substitution they were dealing with.

\medskip

The theorem has the following topological interpretation ( for proofs, see [Qu] or [Ho]):
\medskip

{\bf Corollary 12.} {\sl 

$\imath$) The map $\sigma\colon\xs\to\xs$ is
one--to--one and open.

$\imath\imath$) Every $x\in \xs$ can be written in a unique way 
$$x=\ts^k\bigl(\sigma(y)\bigr) \hbox{ with } y\in\xs\hbox{ and }
0\leq k<|\sigma(y_0)| \ .$$

$\imath\imath\imath$) For every $x\in\xs$, $|\sigma(x_0)|$ is the first return time of
$\sigma(x)$ to $\sigma(\xs)$, i.e. the smallest positive integer $k$ such that
$T^k_{\sigma}\sigma(x)\in\sigma(\xs)$.

$\imath\nu$) The map $\sigma\colon\xs\to\sigma(\xs)$ is an
isomorphism of the system $(\xs,\ts)$ onto the system induced by
$\ts$ on $\sigma(\xs)$. } 

\medskip

The same results hold for $\sigma^n$ for each $n>0$. 
For every letter $a\in A$, we write $[a]=\{x\in \xs\semi x_0=a\}$.
From the  preceding corollary we get immediately:

\medskip

{\bf Corollary 13.} {\sl For every $n>0$,  
$$\P_n=\bigl\{ T^k\bigl( \sigma^n([a])\bigr)\semi a\in A,\;
0\leq k<|\sigma^n(a)|\bigr\}$$
 is a clopen partition of $\xs$.}

\medskip

It is worth noticing that the partition $\P_n$ is a {\sl Kakutani--Rohlin
partition} in the sense of [HPS, section 4]. The {base} of this partition is
$$\bigcup_{a\in A}\sigma^n([a])=\sigma^n(\xs)\ .$$

\medskip

{\bf Proposition 14.} {\sl The sequence of partitions $(\P_n)$ is
nested (see subsection 1.3) i.e.:

$\imath$)  The sequence of bases $(\sigma^n(\xs) ; n\in \N)$ is decreasing.

$\imath\imath$) For every $n$, $\P_{n+1}\succeq \P_n$ as partitions.

Moreover, if the substitution $\sigma$ is proper, then

$\imath\imath\imath$) The intersection of the bases consists in 
only one point (which is the unique fixed point of $\sigma$).

$\imath\nu$) The sequence of partitions spans the topology of
$\xs$.}

\medskip

{\bf Proof.} $\imath$) is obvious. Remember that, for every $x\in\xs$,
$$
\sigma^n(T_{\sigma} x)=T_{\sigma}^{\textstyle |\sigma^n(x_0)|}\sigma^n(x)\ .
$$
Let $a$ be a letter, and $k$ an integer with $0\leq k<|\sigma^{n+1}(a)|$. Set $b_1\ldots
b_m=\sigma(a)$, so that we have 
$$
|\sigma^n(b_1\ldots b_m)|=|\sigma^{n+1}(a)|>k\ .
$$
Let $j\in \{ 0,\ldots m-1 \}$ be the integer defined by
$$
|\sigma^n(b_1\ldots b_j)|\leq k<|\sigma^n(b_1\ldots b_{j+1})| \ .
$$
It follows immediately that 
$$
T^k\sigma^{n+1}([a])\subset T^\ell \sigma^n([b_j])
$$
where $\ell=k-|\sigma^n(b_1\ldots b_j)|$, and $\imath\imath$) is proved.

\medskip

We assume now that $\sigma$ is proper. Then $\imath\imath\imath$) is obvious.
Let us prove  $\imath\nu$):

Let $p,r,\ell$ be as in the definition of a proper substitution (Definition 10).

 Given an integer $m>0$, we claim that, for $n$ large enough,
$x_{[-m,m]}$ is constant on each element of $\P_n$. 

 For $n\geq p$, we write
$R_n=|\sigma^{n-p}(r)|$ and $L_n=|\sigma^{n-p}(\ell)|$. 
Choose $n$ so large that $R_n$ and 
$L_n$ are both greater than $m$.
Fix $a\in A$ and $0\leq k<|\sigma^n(a)|$.  For each $x\in T^k\sigma^n([a])$, there
exists $y\in\xs$ with $y_0=a$ and $x=T^k\sigma^n(y)$.  
The word $\sigma^n(a)\sigma^n(y_1)$ is a prefix of
$\bigl(\sigma^n(y)\bigr)_{[0,\infty)}$, and $\sigma^{n-p}(\ell)$ is a prefix
of $\sigma^n(y_1)$, thus $\sigma^n(a)\sigma^{n-p}(\ell)$ is a prefix of 
$\bigl(\sigma^n(y)\bigr)_{[0,\infty)}$. The same way, 
$\sigma^{n-p}(r)$ is a suffix of $\bigl( \sigma^n(y))_{(-\infty,-1]}$, and
$$
\bigl(
\sigma^n(y)\bigr)_{\textstyle [-R_n,K_n+L_n)}=\sigma^{n-p}(r)\sigma^n(a)
\sigma^{n-p}(\ell)
$$
 where $K_n=|\sigma^n(a)|$, and we get:
$$
x_{[-m,m]}=\Bigl( \sigma^{n-p}(r)\sigma^n(a)\sigma^{n-p}(\ell)\Bigr)_%
{\textstyle [R_n+k-m,R_n+k+m]}
$$
which does not depend on $x$, but only on $a$ and $k$: our claim is proved,
and $\imath\nu$) follows.\cqfd

%% file: section3.tex
{\gros 3. From Bratteli diagrams to substitutions.}

\bigskip

In this section, we prove the first part of Theorem 1. More
precisely, we show here that:
\medskip

{\leftskip 1cm
{ \sl The system associated to a stationary, properly ordered Bratteli diagram is isomorphic, either to the substitution
dynamical system associated to some proper substitution, or to an
odometer with a stationary base.} \par}
\medskip

{\bf 3.1. The substitution read on a stationary ordered Bratteli
diagram.}
\medskip

Let $\B$ be a stationary, properly ordered Bratteli diagram.
Remember that, for each integer $n\geq 0$, $V_n$ denotes the set of vertices
at level $n$. Let us choose a stationary labeling of $V_n$ (for $n>0$) by an
alphabet $A$, i.e. $V_n = \{ V(n,a) ; a\in A \}$ for all $n\in \N$. For $n\geq 1$ and
$a\in A$, $V(n,a)$ is the vertex of label $a$ at level $n$.
 Fix an integer $n>1$. For every letter $a\in A$, consider the
ordered list $(e_1,\ldots ,e_k)$
of edges which range at $V(n,a)$, 
and let $(a_1,\ldots ,a_k)$ be the ordered list of the labels of the sources
of these edges. The map $a\mapsto a_1\ldots a_k$ from $A$ to $A^+$
doesn't depend of $n$. We consider it as a substitution on the alphabet $A$,
and call it {\sl the substitution read on} $\B$.

\medskip
{\bf Lemma 15.} {\sl The substitution $\sigma$ read on $\B$ is primitive and
proper.}

\medskip
{\bf Proof.} As $\B$ is simple, there exists an integer $n>1$ such that, for
every $a,b\in A$, the vertex $V(1,b)$ is connected to the vertex $V(n,a)$; by
definition of $\sigma$, it means that $b$ occurs in $\sigma^{n-1}(a)$: the
substitution $\sigma$ is primitive.

For each $a\in A$, let $i(a)$ be the first letter of $\sigma(a)$: for
each $n>1$, the minimum edge which ranges at  $V(n,a)$ sources at
$V(n-1,i(a))$. Suppose that, for every $n\geq 1$, the range of the map $i^n:A\to
A$ contains at least two letters. Then there exist two distinct sequences 
$(a_n\semi n\geq 1)$ and $(b_n\semi n\geq 1)$ of letters, with $i(a_{n+1})=a_n$
and $i(b_{n+1})=b_n$ for every $n$; by definition of $i$, there exist two
minimal paths $x,y\in\xb$ such that, for every $n\geq 1$, $x$ goes through 
$V(n,a_n)$ and $y$ through $V(n,b_n)$. But $\B$ is properly ordered, and we get a
contradiction. We have showed that, for $p$ large enough, the first letter
of $\sigma^p(a)$ doesn't depend on $a$, and by the same method the last
letter of $\sigma^p(a)$ doesn't depend on $a$: $\sigma$ is proper.
\cqfd

\medskip

{\bf 3.2. A particular case.}

\medskip

Let  $(d_n;n\in \N)$ be a sequence of positive integers. We recall that the inverse limit of the sequence of groups $(\Z/ d_0 d_1\cdots d_n \Z)$, properly topologized, endowed with the addition of $1$ is called the {\it odometer with base $(d_n;n\in \N)$}. If $d=d_n = d_{n+1} = d_{n+2} = \cdots$ for a certain $n$, then we say that the odometer is stationary with stationary base $d$.

\medskip

 We can now prove the result in a particular case:

\medskip

{\bf Proposition 16.} {\sl Let $\B$ be a stationary, properly ordered Bratteli diagram with only simple edges between the top vertex and
the first level, and $\sigma: A\to A^+$ the substitution read on $\B$.

$\imath$) If $\sigma$ is aperiodic, then the system $(\xb,\vb)$ is
isomorphic to the system $(\xs,\ts)$.

$\imath\imath$) If $\sigma$ is periodic, the system  $(\xb,\vb)$ is
isomorphic to an odometer with a stationary base.
}

\medskip

{\bf Proof of $\imath$)}: Consider  the sequence $(\P_n)$ of partitions
associated to $\sigma$ as in 2.4. The ordered Bratteli diagram associated to this
sequence of partitions as in 1.3 is clearly $\B$. By Proposition 14 and Lemma 15, this sequence of partitions satisfies the hypotheses 
$\imath\imath\imath$) an $\imath\nu$) of subsection 1.3, thus
$(\xs,\ts)$ is isomorphic to $(\xb,\vb)$.

This proves $\imath$), but we find it interesting to give an independent proof,
based on an explicit construction of an isomorphism. The notations
introduced here will be used in the proof of $\imath\imath$).

Let $\pi\colon\xb\to A^\Z$ be defined by:

$$
\forall x\in\xb,\;\forall k\in\Z,\,\forall a\in A,\ 
\bigl(\pi(x)\bigr)_k=a\hbox{ if }\vb^k x\hbox{ goes through }V(1,a)\ .
$$
We claim that $\pi$ is an isomorphism between $(\xb,\vb)$ and $(\xs,\ts)$.

The proof of this claim uses some definitions, which will be also used in the
proof of $\imath\imath$). So we don't suppose that $\sigma$ is aperiodic
for the moment.

$\pi$ is clearly continuous, and
$$ \pi\circ\vb=T\circ\pi \eqno{(1)} $$
where $T$ is the
shift on $A^\Z$. Remember that $\xb$ is naturally endowed with a
sequence $(\QQ_n)$ of Kakutani--Rohlin
partitions where $\QQ_n$  is 
defined by the first $n+1$ levels: two paths $x,y\in\xb$ belong to the same
element of $\QQ_n$ if they agree until the level $n+1$. In our case, $\QQ_n$
is easy to describe: Let $K(n,a)$ be
the set of paths $x\in\xb$ such that:

{\parindent=2cm
$\alpha$) $x$ goes through $V(n+1,a)$;

$\beta$) From the top vertex to the level $n+1$, $x$ consists of only minimal edges.
\par}
By induction, the number of paths from the origin to the vertex $V(n+1,a)$
is equal to $|\sigma^n(a)|$ for all $a\in A$ and all $n\geq 0$, and 
 the partition $\QQ_n$ can be written:
$$\QQ_n=\bigl\{ \vb^kK(n,a)\semi a\in A,\;0\leq k<|\sigma^n(a)|\bigr\}\ .$$

 We now define a map $f\colon\xb\to\xb$. Let $x\in\xb$ be a
path;  for $n\geq 1$, let $x(n)$ be  
 the label of the edge between the levels
$n-1$ and $n$ the path $x$ goes through;
 remark that, for the class of Bratteli diagrams we consider,
the sequence $(x(n)\semi n\geq 2)$ completely determines $x$: the edge $x(1)$
of $x$ between the top vertex and the first level is determined by its range,
which is the source of $x(2)$. For all $n>2$,
let $y(n)$ be $x(n-1)$ and $y(2)$ be the label of the minimal element in $E_2$ of the set of edges which range at the source of $x(2)$. The sequence $(y(n)\semi n\geq 1)$ determines a path
$y\in\xs$, and we define $f(x)=y$.

Clearly we have
$$ f(\xmin)=\xmin \ . \eqno{(2)}$$
Suppose that the path $x\in\xb$ goes through $V(1,a)$.  By construction
of the Vershik map we have
$$ f(\vb x)=\vb^{|\sigma(a)|}f(x) \hbox{ and }
f^n(\vb x)=\vb^{\textstyle|\sigma^n(a)|}f(x) \hbox{ for all }n\geq 1\ ;
\hfill\eqno{(3)}$$
moreover, if $\sigma(a)=b_1\ldots b_m$, then for $0\leq j<|\sigma(a)|$ 
the path $\vb^jf(x)$ goes through $V(1,b_{j+1})$; therefore we have

$$ \pi\circ f=\sigma\circ\pi\ .\eqno{(4)}$$
It follows that $\sigma(\pi(\xmin))=\pi(\xmin)$, thus $\pi(\xmin)$ is the
unique fixed point $\omega$ of $\sigma$. By (1), (2), and continuity we get
$\pi(\xb)\subset\xs$, and  $\pi(\xb)=\xs$ by minimality. 

 By definition of $f$, for every letter $a$ we have
$f\bigl(K(0,a)\bigr)=K(1,a)$, and $f^n\bigl(K(0,a)\bigr)=K(n,a)$ for every
$n\geq 1$ by induction. From (1) and (4) we get:
$$ \forall a\in A,\, \forall n\geq 1,\, \forall j\hbox{ with
}0\leq j<|\sigma^n(a)|,\  \pi\bigl(\vb^jK(n,a)\bigr)\subset \ts^j\sigma^n([a])\ .
\eqno{(5)}$$
Suppose now that $\sigma$ is aperiodic.  By (5), $\pi$ maps different
elements of the partition $\QQ_n$ in different elements of the partition
$\P_n$; as the sequence $(\QQ_n)$ of partitions spans the topology of
$\xb$, $\pi$ is one-to-one and our claim is proved. \cqfd

\medskip

{\bf Proof of $\imath\imath$)}: We suppose that $\sigma$ is periodic. Let
$p$ be the smallest period of its fixed point $\omega$, and
$u=\omega_{[0,p)}$, so that $\omega$ is an infinite concatenation of the
word $u$. As $\sigma(\omega)=\omega$, it is also an infinite concatenation
of the word $\sigma(u)$, and it follows that $\sigma(u)=u\ldots u$ ($d$
times) for some integer $d>1$. Thus, for every $a\in A$,
$$ \sum_{b\in A}M(\sigma)_{b,a} |u|_b= |\sigma(u)|_a=d|u|_a$$
where $M(\sigma)$ is the matrix of $\sigma$, defined by
\medskip

\centerline{ for all $a,b\in A$, $M(\sigma)_{b,a}$ is the number of
occurrences of $a$ in $\sigma(b)$.}
\medskip

In other words, the composition vector $(|u|_a\semi a\in A)$ of $u$ is a left
eigenvector of $M(\sigma)$ for the eigenvalue $d$. Thus, for every $n\geq 0$ and
every $a\in A$,
$$
|\sigma^n(u)|_a= d^n|u|_a\hbox{ and } |\sigma^n(u)|=pd^n.
$$
As $\pi$ maps $\xb$ onto $\xs$ and $K(0,a)$ onto $[a]$ for each $a$, the
partition $\QQ_0$ is periodic of period $p$ for $\vb$; it means that every
$\vb$--orbit visits the elements of $\QQ_0$ with period $p$; moreover,
each interval of length $p$ in each $\vb$--orbit visits the set $K(0,a)$ $|u|_a$ times.

For $n\geq 1$,  let $K_n$ be the base of the partition $\QQ_n$, i.e. the union of
the sets $K(n,a)$, and $W_n$ be the transformation induced by $\vb$ on
$K_n$. For each $x\in K(0,a)$, $f^n(x)\in K(n,a)$ and the first
return time of $f^n(x)$ to $K_n$ is $|\sigma^n(a)|$, thus
$$
W_nf^n(x)=\vb^{\textstyle |\sigma^n(a)|}f^n(x)=f^n(\vb x)
$$
by (3), and $f^n$ is an isomorphism of $(\xb,\vb)$ to $(K_n,W_n)$. This
isomorphism maps the partition $\QQ_0$ to the partition
$\RR_n=\{ K(n,a)\semi a\in A\}$ of $K_n$. Therefore, the partition 
$\RR_n$ is periodic for $W_n$ with period $p$,  and every interval of
length $p$ of each $W_n$--orbit visits the set $K(0,a)$ $|u|_a$ times. As
$(K_n,W_n)$ is the system induced by $(\xb,\vb)$ on $K_n$, it
follows that the partition $\QQ_n$ is periodic for $\vb$, with a period
equal to
$$
\sum_{a\in A}|\sigma^n(a)|\,|u|_a=|\sigma^n(u)|=pd^n\ .
$$
Finally, the sequence of partitions $\QQ_n$ is nested, and spans
the topology of $\xb$; $\QQ_n$ is periodic of period $pd^n$: it follows that 
$(\xb,\vb)$ is isomorphic to the odometer with base $(p,d,d,\ldots)$. \cqfd

\medskip

{\bf 3.3. The general case.}

\medskip

We turn now to  the general case, i.e. without the
assumption that the edges from the top vertex to the first   level  are simple.
The result follows immediately from [Fo, Lemma 15] that we recalled in subsection 1.5 (Lemma 9).

%% file: section4.tex
{\gros 4. Return words and derivatives of a sequence}
\bigskip

We present here a modification of the method of  {\sl return
word and derivative sequences}  introduced in [Du1].
This method will be used in the next section to prove that every
substitution dynamical system is isomorphic to the system associated to
some stationary, properly ordered Bratteli diagram; and also in section
7 in the study of factors of substitution dynamical systems. As it
can be interesting in itself, we present it here with more 
generality than we actually need. The proofs which are not given
here can be found in [Du1,Du2]. 

\medskip

{\bf 4.1. Return words.}

\medskip

In all this section, $(X,T)$ is a minimal subshift on the alphabet $A$, and
$x$ a given point of $X$; $x$ is  uniformly recurrent  (see subsection 2.2).

Let $u$ be a suffix of
$x_{(-\infty,-1]}$ and $v$ a prefix of $x_{[0,+\infty)}$, and
we assume that at least one of them is not the empty word.

We define an {\sl occurrence of $u.v$ in $x$} to be an integer $n$ such that
$x_{[n-|u|,n+|v|)}=uv$.
\medskip

{\bf Definition 11.} {\sl A word $w$ on $A$ is a {\rm return word to $u.v$}
in $x$ if there exist two consecutive occurrences $j,k$ of $u.v$ in $x$
such that $w=x_{[j,k)}$. }
\medskip

As $x$ is uniformly recurrent, the
difference between two consecutive occurrences of $u.v$ in $x$ is
bounded, and the set $\rruv$ of return words to $u.v$ is finite. It is
immediate to check that a word $w\in A^+$ is a return word if and only if:

{\parindent=2cm
$\imath$) $uwv\in\L(x)$ (i.e. $uwv$ is a factor of $x$);

$\imath\imath$) $v$ is a prefix of $ wv$ and $u$ is a suffix of
$uw$;

$\imath\imath\imath$) the word $uwv$ contains only two
occurrences of $uv$.
\par}
\medskip
{\bf Remarks and notations.}

{\bf 1)} The statement $\imath\imath$) cannot
be simplified: it is not equivalent to 
{\sl $v$ is a prefix and $u$  a suffix of $w$}.
For example, if  $aaaaa$ is a factor of $x$ then the word $a$
is a return word to $aa.aa$.

{\bf 2)} From this characterization, it follows that the set of return words
doesn't depend on the choice of the point $x$, but only on the subshift
$X$. 

{\bf 3)} If the sequence $x$, or the subshift $(X,T)$ in which we consider 
return words, is not clear from the context, we write $\RR_{u.v}(x)$ or 
$\RR_{u.v}(X)$ instead of $\RR_{u.v}$.

{\bf 4)} To avoid unnecessary heavy notations, we write $\RR_u$ instead of
$\RR_{\emptyset.u}$.
\medskip

{\bf Lemma 17.} {\sl The set $\RR_{u.v}$ is a code. Moreover, it is a circular
code (Definition 9)}.
\medskip

{\bf Proof.} Let $m$ be a word, and assume that $m$ can be written an a
concatenation $m=w_1\ldots w_k$ of return words. 
From the property $\imath\imath$) of return words, it follows by
induction on $j$ that $u$ is a suffix of $uw_1\ldots w_j$ for $1\leq j\leq k$.
By backwards induction, $v$ is a prefix of $w_j\ldots w_kv$ for
 $1\leq j\leq k$. Thus $u.v$ has at least $k+1$ occurrences in $umv=uw_1\ldots
w_kv$, namely between $u$ and $w_1$; between $w_j$ and $w_{j+1}$ for
$1\leq j<k$; and between $w_k$ and $v$. By the property 
$\imath\imath\imath$) of return words, $u.v$ cannot have other
occurrences in $umv$. The decomposition $m=w_1\ldots w_k$ of $m$ is
obtained in cutting the word $umv$ at each occurrence of $u.v$; this
decomposition is therefore unique, and $\RR_{u.v}$ is a code. 
Moreover, $u$ is a suffix of
$w_1\ldots w_k$  if  $|w_1\ldots w_k|\geq |u|$; 
and $v$ a prefix of $w_1\ldots w_k$   if $|w_1\ldots w_k|\geq |v|$. It follows
that if a  sequence $y\in A^\Z$ can be decomposed as a concatenation
of return words, the occurrences of $u.v$ in $y$ are exactly the cutting
points of this decomposition. From that we deduce easily that $\RR_{u.v}$
is a circular code:

Let $w_1,\ldots,w_j,w,w'_1,\ldots,w'_k\in \RR_{u.v}$, $s,t\in A^*$
be as in the definition of a circular code, and  $y$ be the periodic point
$$y=(w_1\ldots w_j)^\infty\,\mid\, (w_1\ldots w_j)^\infty \ .$$
By the discussion above, $0$ is an occurrence of $u.v$ in $y$. But
$$y=(ww'_1\ldots w'_k)^\infty\,t\,\mid\,s\,(w'_1\ldots w'_kw)^\infty$$
and $ts=w\in\RR_{u.v}$; by the discussion above again, $0$ is not an
occurrence of $u.v$ in $y$, except if $t$ is the empty word.  \cqfd

\medskip

{\bf 4.2. Derivatives of a sequence.}

\medskip

It will be convenient to label the return words. Put
$$\ruv = \bigl\{ 1,\ldots,\card(\rruv)\bigr\}$$
and let $\phi_{u.v}: R_{u.v}\to\RR_{u.v}$ be the bijection defined as
follows: let $\RR_{u.v}$ be ordered according to the rank of first occurrence
in $x_{[0,+\infty)}$, and  $\phi_{u.v}(k)$ defined to be the 
$k^{{\hbox{\petit th}}}$ element of $\RR_{u.v}$ for this order. When $u$ is the empty-word we set $\phi_{u.v} = \phi_v$.

We consider $\ruv$ as an alphabet, and $\phiuv$ as a map from
$\ruv$ to $A^+$. The last lemma can be stated as follows:
\medskip

{\bf Corollary 18. } {\sl $\phi_{u.v}:R_{u.v}^+\to A^+$ and 
$\phi_{u.v}:R_{u.v}^\Z\to A^\Z$ are one to one.}
\medskip

The sequence $x$ itself is by construction a  concatenation of return
words, thus we can define:
\medskip

{\bf Definition 12.} {\sl The $u.v$--{\rm derivative} of $x$ is the unique
sequence $\D_{u.v}(x)$ on the alphabet $R_{u.v}$ such that:}
$$\phi_{u.v}\bigl( \D_{u.v}(x)\bigl)=x \ .$$

\medskip

{\bf 4.3. Topological interpretation.}

\medskip

The notion
of derivative sequence is the combinatorial analogue to the notion of
induced system, as we explain now.

Let $u$, $v$ be as above. The cylinder sets $[u.wv]$ for $w\in\RR_{u.v}$
are obviously pairwise disjoint; they are included in the cylinder set
$[u.v]$ by the property $\imath\imath$) of return words. Let $y\in [u.v]$,
and $n$ be the smallest positive occurrence of $u.v$ in $y$; then
$w=y_{[0,n)}$ is a return word, and $y\in [u.wv]$. Thus 
$\{ [u.wv]\semi w\in\RR_{u.v} \}$ is a partition of $[u.v]$. Moreover, if
$w\in\RR_{u.v}$ and  $y\in [u.wv]$, the first return time of $y$ to
$[u.v]$ is $|w|$ by the property $\imath\imath\imath$) of return words. It
follows that
$$
\QQ=\bigl\{ T^j[u.wv] \semi w\in\RR_{u.v}\hbox{ and }0\leq j<|w| \bigl\}
$$
is a Kakutani--Rohlin partition of $X$, with base $[u.v]$.
Using the bijection $\phi_{u.v}$, this partition can also be written:
$$
\QQ=\bigl\{ T^j[u.\phi_{u.v}(k)v]\semi k\in R_{u.v}\hbox{ and }
0\leq j<|\phi_{u.v}(k)|\bigr\}\ .
$$
Let $S$ be the shift on $\ruv^\Z$ and $Y$  the subshift spanned by
$\duv(x)$. 
\medskip

{\bf Lemma 19.} {\sl $\phiuv$ is an isomorphism of $(Y,S)$ onto the system
induced by $(X,T)$ on the cylinder set $[u.v]$.}

\medskip
{\bf Proof.} We know that $\phiuv\colon R_{u.v}^\Z\to
A^\Z$ is one to one.
As
$$
\phiuv(\duv(x))=x\hbox{ and, for all } y\in \ruv^\Z, \ 
\phiuv(Sy)=T^{\textstyle|\phiuv(y_0)|}\,\phiuv(y)
$$
we get $\phiuv(S^n\duv(x))\in X$ for all $n$, thus $\phiuv(Y)\subset X$. By
definition of $\phiuv$, we have $\phiuv(Y)\subset [u.v]$.

 Let $z\in [u.v]$. There exists a  sequence $(n_i)$ of integers such that
$T^{n_i}x\to z$; as $[u.v]$ is open in $X$, for $i$ large enough $T^{n_i}x\in
[u.v]$, and $n_i$ is an occurrence of $u.v$ in $x$; it follows that
$T^{n_i}x=\phiuv(S^{k_i}\duv(x))$ for some $k_i$, and
$T^{n_i}x\in\phiuv(Y)$; finally we get $z\in\phiuv(Y)$ and   
$\phiuv(Y)=[u.v]$.

Let $y\in Y$ and $z=\phiuv(y)$; the first return time of $z$ to $[u.v]$ is
$n=|\phiuv(y_0)|$; thus the image of $z$ by the first return time
transformation is $T^nz=\phiuv(Sy)$, and the lemma is proved. \cqfd
\medskip

{\bf 4.4. Return words and Bratteli diagrams.}
\medskip

We consider now return words for longer and longer prefixes of
$x_{[0,+\infty)}$ and suffixes of $x_{(-\infty,-1]}$. To avoid unnecessary
heavy notations, we write, for all $n\geq 1$:
$$
\RR_n=\RR_{\textstyle x_{[-n,-1]}.x_{[0,n)}}\semi
R_n=R_{\textstyle x_{[-n,-1]}.x_{[0,n)}}\semi 
\phi_n=\phi_{\textstyle x_{[-n,-1]}.x_{[0,n)}} \ .
$$
We get a sequence of Kakutani--Rohlin partitions $\P_n$ where
$$
\P_n  = \Bigl\{  T^j\bigl[  x_{[-n,-1]}.wx_{[0,n)}\bigr]  \semi
w\in\RR_n\hbox{ and } 0\leq j<|w| \Bigl\} 
$$    
$$
 =  \Bigl\{ T^j B_{n,k} \semi
 k\in R_n \hbox{ and } 0\leq j<|\phi_n(k)| \Bigr\} 
$$

where  $B_{n,k}  =\bigl[ x_{ [-n,-1] }.\phi_n(k)x_{[0,n)} \bigr] 
\hbox{ for }k\in R_n\ .$
 The base of $\P_n$ is
 $B_n=\bigl[ x_{[-n,-1]}.x_{[0,n)}\bigr]$.

\medskip
Fix $n\geq 1$. As $x_{[-n,-1]}$ is a suffix of $x_{[-n-1,-1]}$, and 
$x_{[0,n)}$  a prefix of $x_{[0,n+1)}$, each word belonging to $\RR_{n+1}$
can be written a concatenation of words belonging to $\RR_n$, and it can
be done in a unique way. It follows that there exists a unique map
 $\sigma_n: R_{n+1}\to R_n^+$ with:
$$\phi_n\circ \sigma_n=\phi_{n+1}\ .$$

 Clearly $B_{n+1}\subset B_n$. Let $k$ belong to $R_{n+1}$, $j_1\ldots j_m=\sigma_n(k)$ and $y\in B_{n+1,k}$. By definition, $y\in B_{n,j_1}$.
The orbit of $y$ climbs the  $j_1^\th$ tower of $\P_n$, then goes to 
$B_{n,j_2}$, climbs the  $j_2^\th$ tower of $\P_n$ \dots\ goes to
$B_{n,j_m}$, climbs the  $j_m^\th$ tower of $\P_n$ and comes back to
$B_{n+1}$. The sequence of elements of $\P_n$ the orbit of $y$ visits
before it returns to $B_{n+1}$ doesn't depend on the point $y\in B_{n+1,k}$,
thus the partition $\P_{n+1}$ is finer than $\P_n$. The sequence $(\P_n)$
of partitions is nested in the sense of subsection 1.3. 

Moreover, $\bigcap_n B_n$ consists exactly of the point $x$, and the
sequence $(\P_n)$ spans the topology of $X$ because $y_{[-n,n)}$ is
constant on each element of $\P_n$.  Thus this sequence satisfy all the
hypotheses $\imath$), $\imath\imath$), $\imath\imath\imath$) and $\imath\nu)$ of subsection 1.3: it can be used to construct a
properly ordered Bratteli diagram $\B$ with $(\xb,\vb)$ isomorphic to
$(X,T)$, or to compute the dimension group of $(X,T)$. It may be interesting to
remark that the incidence matrices of this diagram --- which are used to
compute the dimension group --- are the matrices $M(\sigma_n)$
associated to the maps $\sigma_n$ and defined by
$$
\hbox{For }  i\in R_{n+1}  \hbox{ and }   j\in R_n\, , \ 
\Bigl( M( \sigma_n )\Bigr)_{i,j} =\hbox{ number of occurrences of } j
\hbox{ in } \sigma_n(i)\ .
$$

%% file: section5.tex
{\gros 5. From substitutions to Bratteli diagrams.}
\bigskip

{\bf 5.1.} In this section, we prove the second part of Theorem 1 in
the form of the following proposition:

\medskip

{\bf Proposition 20.}
{\sl The system associated to any
(primitive, aperiodic) substitution is isomorphic to the system associated
to some stationary, properly ordered Bratteli diagram.}
\medskip

Suppose first that the  (primitive, aperiodic)
substitution $\sigma$ is  proper.  We can clearly define a stationary
ordered Bratteli diagram $\B$, with simple edges between the top vertex and
the first level, such that the substitution read on $\B$ in the sense of
Section 3 is exactly $\sigma$. This diagram is simple because $\sigma$ is
primitive, and properly ordered because $\sigma$ is proper. As $\sigma$ is aperiodic, Proposition 16 of Section 3 says
 that $(\xs,\ts)$ is isomorphic to $(\xb,\vb)$, and our result is proved. 
We consider now the general case, without assuming that the substitution
is proper.

\medskip

{\bf Remark.} The proposition  is partially proved in [Fo]. Our approach
differs from Forrest's in that
we give an explicit (algorithmic) construction of  a properly ordered Bratteli
diagram, and also of a proper substitution, defining a system isomorphic
to $(\xs,\ts)$; in particular, it
leads to an explicit method to compute completely (i.e.: with the order
and the order unit) the dimension group of any 
substitution dynamical system:
 we don't make any assumption on the given substitution ---
except of course that it is primitive and aperiodic.

\medskip

{\bf 5.2. A derivative of a fixed point of a substitution.}

\medskip

Let $\sigma$ be a (primitive, aperiodic) substitution on the alphabet $A$,
and $x$ one of its fixed points. We write $r=x_{-1}$ and $\ell=x_0$, so
that $r$ is the last letter of $\sigma(r)$ and $\ell$ the first letter of
$\sigma(\ell)$. 

We use here freely the notations introduced in
Section 4. But, as we consider here return words to $r.\ell$ only,  we
write $\RR$, $R$, $\phi$, $\D(x)$ instead of $\RR_{r.\ell}$, $R_{r.\ell}$,
$\phi_{r.\ell}$ and $\D_{r.\ell}(x)$ respectively.

We define now a substitution $\tau$ on the alphabet $R$.

 Let $j\in R$, and  $w=w_1\ldots w_k=\phi(j)\in \RR$.
By the characterization  of return words (see subsection 4.1),
 we have  $rw\ell\in\L(\sigma)$,
$w_1=\ell$ and $w_k=r$.
As $\sigma(x)=x$, the word 
$\sigma(rw\ell)=\sigma(r)\sigma(w)\sigma(\ell)$
belongs also to $\L(\sigma)$.
But 
$r$ is the last letter of $\sigma(r)$, and $\ell$ the first of
$\sigma(\ell)$; it follows that $r\sigma(w)\ell\in\L(\sigma)$; and also
that the first letter of $\sigma(w)$ is $\ell$ and its last letter  $r$:
$r\ell$ is a prefix and also a suffix of $r\sigma(w)\ell$. Therefore, the
word $\sigma(w)$ appears in $x$ between two occurrences of $r.\ell$, thus
it is a concatenation of return words, i.e. belongs to $\phi(R^+)$; as $\phi:
R\to A^+$ is a  code (Lemma 17), there
exists a unique word $u\in R^+$ such that $\sigma(w)=\phi(u)$. 

We define $\tau(j)=u$;
 $\tau$ is a substitution on the alphabet $R$, characterized by
$$\phi\circ\tau=\sigma\circ\phi\ .$$
It follows that $\phi\circ\tau^n=\sigma^n\circ\phi$ for each $n\geq 0$.

\medskip

{\bf Lemma 21.} {\sl  $\tau$ is proper, primitive, aperiodic, and $\D(x)$ is its
fixed point.}

\medskip

{\bf Proof.} Recall that $\phi(1)$ is the first element in the decomposition
of $x_{[0,+\infty)}$ in return  words, i.e. that $\phi(1)\ell$ is a prefix of 
$x_{[0,+\infty)}$. Let $n$ be so large that $|\sigma^n(\ell)|>|\phi(1)|$. As
$\sigma^n(\ell)$ and $\phi(1)\ell$ are both prefixes of $x_{[0,+\infty)}$, 
$\phi(1)\ell$ is a prefix of $\sigma^n(\ell)$. Let $j\in R$, and $w=\phi(j)$.
As $\ell$ is the first letter of $w$, $\sigma^n(\ell)$ is a prefix of
$\sigma^n(w)$, thus $\phi(1)\ell$ also. It follows that $\phi(1)$ is the
first element in the decomposition of $\sigma^n(w)=\phi(\tau^n(j))$ in a
concatenation of return words, i.e. that $1$ is the first letter of $\tau^n(j)$.

Let $m=\bigl(\D(x)\bigr)_{-1}$: $r\phi(m)$ is a suffix of
$x_{(-\infty,-1]}$, and the same argument shows that, for every $n$ large
enough and every $j\in R$, $m$ is the last letter of $\tau^n(j)$: 
$\tau$ is proper.

Let $k>0$ be an occurrence of $r.\ell$ so large that every return word
$w\in\RR$ appears in the decomposition of $x_{[0,k)}$, i.e. that every $j\in
R$ occurs in the word $u\in R^+$ defined by $\phi(u)=x_{[0,k)}$;
let $n$ be so large that $|\sigma^n(\ell)|>k$. Let $i,j\in R$. As above,
$x_{[0,k)}\ell$ is a prefix of $\sigma^n(\ell)$, which is a prefix of
$\sigma^n(\phi(i))=\phi(\tau^n(i))$. Thus $u$ is a prefix of $\tau^n(i)$, and
$j$ occurs in $\tau^n(i)$: $\tau$ is primitive.
Moreover
$$\phi(\tau(\D(x)))=\sigma(\phi(\D(x)))=\sigma(x)=x=\phi(\D(x))$$
thus $\tau(\D(x))=D(x)$ by the unique decomposition property, and $\D(x)$
is the unique fixed point of $\tau$. As $\phi(\D(x))=x$ is not periodic,
$\D(x)$ is not periodic and the substitution $\tau$ is aperiodic. \cqfd

\medskip

{\bf 5.3. Proof of the Proposition 20.}

\medskip

The proof of the announced result is almost immediate from this point: 
Let $\A$ be the stationary, properly ordered Bratteli diagram
naturally associated to $\tau$: it has  simple edges between the top level and
the  first level, and $\tau$ is the substitution read on it. 

Let $\B$ be the stationary, properly ordered diagram which is
identical to $\A$ except that, for each $j\in R$, the top level is joined by
$|\phi(j)|$ edges to the vertex $V(1,j)$. We claim that 
$(\xs,\ts)$ is isomorphic to $(\xb,\vb)$.

In Section 3 we have constructed an isomorphism $\pi^{-1}$ from
$(\xtau,\ttau)$ to $(X_\A,V_\A)$. For every $j\in R$, this isomorphism maps
the cylinder set $[j]$ of $\xtau$ to the subset of $X_A$ consisting of
the paths which go through the vertex $V(1,j)$. Thus the systems
$(\xs,\ts)$ and $(\xb,\vb)$ are constructed in putting Kakutani-Rohlin
towers of the same height over  subsets of $\xtau$ and $X_\A$
respectively, which correspond via this isomorphism: hence these systems are
isomorphic.  \cqfd
\medskip

{\bf 5.4. The dimension group of a substitution dynamical system.}
\medskip

The construction above allows us to compute the dimension group of any
substitution dynamical system, including the order and the order unit. 

When $\sigma$ is a substitution on an alphabet $A$, we let $M(\sigma)$ denote its matrix; it is the $|A|\times |A|$ matrix, with rows and
columns indexed by $A$, defined by:
\medskip

\centerline{For $a,b\in A$, $M(\sigma)_{a,b}$ is the number of occurrences
of $b$ in $\sigma(a)$.}
\medskip

Let us recall the definition of the dimension group $(K,K_+,\one)$ of a
$d\times d$ matrix $M$ with nonnegative integer entries. Define:
$$
G = \bigl\{ \um\in\Q^d\semi \exists n>0\hbox{ such that }M^n\um\in
\Z^d\bigr\}
$$
$$
G_+ = \bigl\{ \um\in\Q^d\semi \exists n>0\hbox{ such that }M^n\um\in
\Z_+^d\bigr\} 
$$
$$
H = \bigl\{ \um\in\Q^d\semi \exists n>0\hbox{ such that
}M^n\um = \overrightarrow 0 \bigr\} 
$$
$$
K = G/H \semi K_+ \hbox{is the image of }G_+\hbox{ in } K 
$$
and $\one$ is the image in $K$ of the vector $\overrightarrow
1=(1,\ldots,1)$.

\medskip

{\bf Theorem 22. } {\sl Let $\sigma$ be a (primitive, aperiodic)
substitution on the alphabet $A$.

$\imath$) If $\sigma$ is proper, then $K^0(\xs,\ts)$ is
isomorphic to the dimension group of the matrix of $\sigma$ as ordered
groups with order units.

$\imath\imath$) In the general case, let $\tau$ be the substitution on the
alphabet $R$, and $\phi\colon R\to A^+$ the morphism defined in 5.2. Then
$K^0(\xs,\ts)$ is isomorphic as ordered group to the dimension group of the
matrix of $\tau$; through this isomorphism, the order unit of $K^0(\xs,\ts)$ 
corresponds to the image of the vector $\bigl(|\phi(j)|\semi j\in R\bigr)$ in
the dimension group of this matrix.}

\medskip

{\bf Proof.} If $\sigma$ is proper, then $(\xs,\ts)$ is
isomorphic to $(\xb,\vb)$ where $\B$ is a   stationary, properly
ordered Bratteli diagram (see 5.1.) with simple edges between the top vertex
and the first level, and $\sigma$ is the substitution read on $\B$.  The
incidence matrix of $\B$ is $M(\sigma)$ at each level, and the result
follows immediately from the computation of the dimension group of the
system associated to a Bratteli diagram (Theorem 6). 

In the general case, let $\B$ be the Bratteli diagram constructed in 5.3. 
At each level, its incidence matrix is the matrix of $\tau$;
it has $|\phi(j)|$ edges between the top vertex and the vertex of label $j$ in the
first level, and the result follows. \cqfd

\medskip

{\bf 5.5. An  example.}

\medskip
 
The sets  $\RR$ and $R$,  the map $\phi$ and the substitution
$\tau$ on $R$ can be computed in an explicit way from the substitution
$\sigma$, as shown in the following example.

Let $\sigma$ be the substitution defined on the alphabet $A=\{a,b\}$ by:
$$\sigma(a)= aba\ ;\ \sigma(b)=baab\ .$$
This substitution is clearly primitive, and it can be checked that it is
aperiodic.

Let $x$ be the fixed point of $\sigma$ such that $x_{-1}=b$ and
$x_0=a$: 
$$x = \ldots\; abaababaab\mid ababaababa\;\ldots$$ 
 where the
vertical bar separates $x_{(-\infty,-1]}$ of $x_{[0,+\infty)}$ as usual. As
$ba$ is a factor of $\sigma(a)$, $x$ is admissible.

 As $baba$ is a factor of $x$, $w_1=ab$ is a return word to  $b.a$; $w_1$ is
in fact the first return word in the decomposition of $x_{[0,+\infty)}$, thus
we set $\phi(1)=w_1$. We have:
$$
\sigma(w_1)= ababaab=w_1w_1w_2
$$
where $w_2=aab$ is necessarily a return word; we set $\phi(2)=w_2$, and
we get $\tau(1)= 112$. We have:
$$
\sigma(w_2)= abaababaab=w_1w_2w_1w_2
$$ 
thus $\tau(2)=1212$. As the substitution $\tau$ is primitive, the alphabet
$R$ cannot contain any other letter,  thus $\RR$ cannot contain any other
return word. We get:
$$
R=\{1,2\}\semi \RR=\{ab,aab\}\semi \phi(1)=ab\semi\phi(2)=aab\semi
\tau(1)=112\semi \tau(2)=1212$$

We can check that, as announced, the substitution $\tau$ is primitive,
aperiodic and proper; its fixed point is
$$
y= \ldots\; 12121121212\mid 1121121212\;\ldots
$$
which is equal to $\D(x)$, and thus satisfies $\phi(y)=x$. 

%{\bf Draw here the  Bratteli diagrams $\A$ and $\B$ corresponding to the
%example.}

Let us compute the dimension group of the substitution dynamical system
$(\xs,\ts)$; we use the notations of 5.4..

Here the matrix $M(\tau)$ is 
$$M(\tau)=\pmatrix{ 2 & 1\cr 2 & 2\cr}\ .$$
The determinant of $M(\tau)$ is not $0$, thus $H=\{0\}$ and $
K^0(\xs,\ts)=G$.

As the determinant of $M(\tau)$ is $2$, for every $\um\in G$ there exists
an integer $k>0$ with $2^k\um\in\Z^2$. Moreover, it is easily checked
that, for every integer $k>0$ there exists an integer $n>0$ such that all the
coefficients of $M(\tau)^n$ are multiples of $2^k$. It follows that
$$K^0(\xs,\ts)=
G=\bigl\{ (2^{-k}a,2^{-k}b)\semi a,b\in\Z,\;k\geq 0 \bigr\}\subset \Q^2\ .$$

With $\rho=2+\sqrt 2$ (the Perron-Frobenius eigenvalue) we have
$$
2\sqrt 2 \rho^{-n}M(\tau)^n\to \pmatrix{\sqrt 2 & 1\cr 2 & \sqrt 2}
\hbox{ as }n\to+\infty\ .
$$
It follows that
$$
K^0_+(\xs,\ts)=G_+=\bigl\{ (\alpha,\beta)\in K^0(\xs,\ts)\semi
\alpha\sqrt 2+\beta \geq 0\bigr\}\ .
$$
The order unit is $\one=(2,3)$.

%% file: section6.tex
{\gros 6. Construction of a proper substitution.}

\bigskip

We use here the notations of the preceding section.

By Lemma 9, there exists a stationary, properly ordered Bratteli diagram $\CC$ with simple edges between
the origin and the first level, which is order equivalent to $\B$,
thus defines an isomorphic system. 
From Proposition 16 we know that 
$(\xs,\ts)$ is isomorphic to $(\xz,T_{\zeta})$, where $\zeta$ is the
substitution read on $\CC$; this substitution is proper by Lemma 15. In Figure 4 (Section 1) this is illustrated. In fact the leftmost diagram corresponds to $\B$ (and $\tau$), while the rightmost diagram corresponds to $\CC$ (and $\zeta$).

It is perhaps interesting to give an explicit and direct construction of such a $\zeta$, without the use of Bratteli diagrams. This will provide an
alternative proof of Proposition 20, in
the formalism of substitutions. Moreover, this construction is
perhaps more algorithmic than the one arising from Lemma 9 proved by A. Forrest. The work is done by the following proposition, which is a
modification of an unpublished result of G. Rauzy.
\medskip 

{\bf Proposition 23.} {\sl Let $y$ be an admissible fixed
point of a primitive substitution $\tau$ on the alphabet $R$, $A$ an
alphabet, $\phi\colon R\to A^+$ a one-to-one map, $x=\phi(y)$, and $(X,T)$ the
subshift spanned by $x$. 

There exist a primitive substitution $\zeta$ on an alphabet $B$, an
admissible fixed point $z$ of $\zeta$ and a map $\theta\colon B\to
A$ such that:

$\imath$) $\theta(z)=x$.

$\imath\imath$) If $\phi(R)$ is a circular code, then $\theta$ is an
isomorphism from $(\xz,T_{\zeta})$ to $(X,T)$.

$\imath\imath\imath$) If $\tau$ is proper, then $\zeta$ is proper.
}

From $\imath$) it follows that $\zeta$ is aperiodic whenever the
sequence $x$ is not periodic.

\medskip

{\bf Proof.} The proof below is very simple, but the notations are, in an
inescapable way, a bit heavy. 
 Substituting a power of $\tau$ for $\tau$ if needed, we can assume
that $|\tau(j)|\geq |\phi(j)|$ for all $j\in R$. For all $j$ in $R$, we
write  $m_j=|\tau(j)|$ and $n_j=|\phi(j)|$. We define:

an alphabet $B$ by  $B =\bigl\{  (j,p)\semi j\in R,\;
1\leq p\leq n_j \bigr\} $,
 
a map $\theta\colon B\to A$  by  $\theta(j,p)= \bigl(\phi(j)\bigr)_p $ ,

a map $\psi\colon R\to B^+ $  by  $\psi(j)=(j,1)(j,2)\ldots (j_,n_j)$.

Clearly $\theta\circ\psi=\phi$.
                                              
We define a substitution $\zeta$ on $B$ by: 
$$
\hbox{For $j$ in $R$
and $1\leq p\leq n_j$, }\zeta(j,p) =\left\{  \matrix{
\psi\Bigl(\bigl(\tau(j)\bigr)_p\Bigr) & \hbox{ if }&1\leq p<n_j\cr
\psi\Bigl(\bigl(\tau(j)\bigr)_{\textstyle[n_j,m_j]}\Bigr) & \hbox{ if
}&p=n_j\cr
 }\right.  
$$

So that, for every $j\in R$, $\zeta\bigl(\psi(j)\bigr) =
\zeta(j,1)\ldots\zeta(j,n_j)=\psi\bigl(\tau(j)\bigr)$, i.e. 
$$
\zeta\circ\psi=\psi\circ\tau\eqno{(2)}
$$
 and it follows that  $$\zeta^n\circ\psi=\psi\circ\tau^n \hbox{ for
all }n\geq 0\ .$$ We claim that $\zeta$ is primitive. Let $n$ be an
integer such that $b$ occurs in $\tau^n(a)$ for all $a,b\in R$. Let
$(j,p)$ and $(k,q)$ belong to $B$. By construction, $\zeta(j,p)$
contains $\psi\bigl(\tau(j)_p\bigr)$ as a factor, thus
$\zeta^{n+1}(j,p)$ contains  $\zeta^n\Bigl(
\psi\bigl(\tau(j)_p\bigr)\Bigr)=
\psi\Bigl(\tau^n\bigl(\tau(j)_p\bigr)\Bigr)$ as a factor. By the
choice of $n$, $k$ occurs in $\tau^n\bigl(\tau(j)_p\bigr)$, thus
$\psi(k)$ is a factor of
$\psi\Bigl(\tau^n\bigl(\tau(j)_p\bigr)\Bigr)$, and also of 
$\zeta^{n+1}(j,p)$. As $(k,q)$ is a letter of $\psi(k)$, $(k,q)$
occurs in  $\zeta^{n+1}(j,p)$ and our claim is proved.

Let $z=\psi(y)$. From (2) we get $\zeta(z)=\psi(\tau(y))=\psi(y)=z$,
and $z$ is a fixed point of $\zeta$. By construction, $z$ is
uniformly recurrent, thus it as an admissible fixed point of
$\zeta$. Moreover, $\theta(z)=\theta(\psi(y))=\phi(y)=x$, and
$\imath$) is proved.

{\bf Proof of $\imath\imath$)}. As $\theta$ commutes with the shift
and maps $z$ to $x$, it maps $\xz$ onto $X$. It remains to prove that
$\theta\colon \xz\to X$ is one-to-one. Let $\alpha\in X$. By
definition of $X$, there exist $\gamma\in \xtau$ and an integer $p$,
with $0\leq p<|\phi(\gamma_0)|$, such that $\alpha=T^p\phi(\gamma)$. Let
$\beta$ be an element of $\xz$ with
$\theta(\beta)=\alpha$. By definition of $\psi$, there exist some
$\delta\in \xtau$ and some integer $q$, with $0\leq q<|\psi(\delta_0)|$,
such that $\beta=T_{\zeta}^q\psi(\delta)$; it follows that
$T^q\phi(\delta)=\theta(\beta)=\alpha=T^p\phi(\gamma)$; as
$0\leq q<|\psi(\delta)|=|\phi(\delta)|$ by construction of $\psi$, and
$\phi(R)$ is a circular code, it follows that $\delta=\gamma$ and $q=p$,
thus $\beta=T_{\zeta}^p\psi(\gamma)$: $\beta$ is uniquely determined by
$\alpha$, and $\theta$ is one-to-one.

{\bf Proof of $\imath\imath\imath$)}. Let $\ell\in R$ be the letter
such that $\ell$ is the first letter of $\tau(k)$ for every $k\in
R$. Let $(j,p)\in B$, and $k=\tau(j)_p$. By definition of $\zeta$,
the first letter of $\zeta(j,p)$ is $(k,1)$, and the first letter of
$\zeta^2(j,p)$ is the first letter of $\zeta(k,1)$, i.e. $(\ell,1)$. By
the same method, if $r$ is the last letter of $\tau(k)$ for every
$k\in R$, then the last letter of $\zeta^2(j,p)$ is $(r,n_r)$ for
every $(j,p)\in B$. \cqfd

\medskip

{\bf The example of subsection 5.5 continues.}

\medskip

As $|\tau(j)|\geq |\phi(j)|$ for $j=1,2$, we need not to consider a
higher power for $\tau$. We have $n_1=2$, $m_1=3$, $n_2=3$ and
$m_2=4$.

The alphabet $B$ is: 
$$
B=\bigl\{
(1,1),\;(1,2),\;(2,1),\;(2,2),\;(2,3)\bigr\}
$$ 
The maps $\theta\colon B\to A$ and $\psi\colon R\to B^+$ are given by: $$
\displaylines{
\theta(1,1)=a\semi \theta(1,2)=b\semi\theta(2,1)=a \semi
\theta(2,2)=a\semi \theta(2,3)=b\cr
\psi(1)=(1,1)(1,2)\semi\psi(2)=(2,1)(2,2)(2,3)\ .\cr }
$$ 
We compute now the substitution $\zeta$ on $B$: 
$$
\matrix{ \zeta(1,1)\hfill&=& \psi\bigl(\tau(1)_1\bigr)\hfill &=&
\psi(1)\hfill&=&(1,1)(1,2)\hfill\cr
\zeta(1,2)\hfill&=& \psi\bigl(\tau(1)_{[2,3]}
\bigr)\hfill &=& \psi(12)\hfill&=&(1,1)(1,2)(2,1)(2,2)(2,3)\hfill\cr
\zeta(2,1)\hfill&=&
\psi\bigl(\tau(2)_1\bigr) \hfill&=&
\psi(1)\hfill&=&(1,1)(1,2)\hfill\cr
\zeta(2,2)\hfill&=&
\psi\bigl(\tau(2)_2\bigr)\hfill &=& \psi(2)&=&(2,1)(2,2)(2,3)
\hfill\cr
\zeta(2,3)\hfill&=& \psi\bigl(\tau(2)_{[3,4]}\bigr)\hfill &=&
\psi(12)\hfill&=&(1,1)(1,2)(2,1)(2,2)(2,3)\hfill\cr }
$$

To get easier notations, we put: $$\alpha=(1,1)\semi
\beta=(1,2)\semi \gamma=(2,1)\semi \delta=(2,2) \semi
\epsilon=(2,3)$$ and $B=\{ \alpha,\beta,\gamma,\delta,\epsilon\}$.
The substitution $\zeta$ and the map $\theta\colon B\to A$ can be
written: 
$$
\displaylines{ \zeta(\alpha)=\alpha\beta\semi
\zeta(\beta)=\alpha\beta\gamma\delta\epsilon\semi \zeta(\gamma)=
\alpha\beta\semi \zeta(\delta)=\gamma\delta\epsilon\semi
\zeta(\epsilon)=\alpha\beta\gamma\delta\epsilon\ .\cr
\theta(\alpha)=a\semi \theta(\beta)=b\semi\theta(\gamma)=a \semi
\theta(\delta)=a\semi \theta(\epsilon)=b\ .\cr }
$$
$\zeta$ is proper: for every $m\in B$,  $\alpha$ is the first and $\epsilon$
the last letter of $\zeta^2(m)$. The fixed point of $\zeta$ is:
$$
z=\ldots\;
\alpha\beta\gamma\delta\epsilon\alpha\beta\gamma\delta\epsilon
\alpha\beta\gamma\delta\epsilon\mid
\alpha\beta\alpha\beta\gamma\delta\epsilon\alpha\beta
\alpha\beta\gamma\delta\epsilon\;\ldots
$$ 
We have $\theta(z)=x$; $\theta$ is an isomorphism from $(\xz,T_{\zeta})$ to $(\xs,\ts)$, and hence $(X_{\sigma} , T_{\sigma} )$ is isomorphic to the Bratteli-Vershik system associated to the stationary, properly ordered Bratteli diagram naturally derived from $\zeta$.

%% file: section7.tex
{\gros  7. Cantor factors of substitution subshifts}
 
\bigskip

In this section we prove Theorem 3 and Theorem 4.

\medskip

In subsection 7.1 we introduce the notion of linearly recurrent 
subshifts and prove some properties of such subshifts. In part 7.2
we prove Theorem 3 and we show that a minimal Cantor system
which is a factor of a linearly recurrent subshift is either
isomorphic to a subshift or to an odometer. This last result is the
core of the proof of Theorem 4, given in part 7.3.

\medskip

{\bf 7.1. Linearly recurrent sequences and subshifts.}

\medskip

{\bf Definition 13.} {  \sl We say that a sequence $x$ on a
finite alphabet is {\rm linearly recurrent (with constant $K\in  
\N$)} if it is recurrent and if, for every factor
$u$ of $x$, the difference between two successive occurrences of $u$
in $x$ is less than $K|u|$.}

\medskip

{\bf Theorem 24.} {\sl  Let $x$ be a linearly recurrent
aperiodic sequence with constant $K$. Then:                  

$\imath$) The number of distinct factors of length $n$ of
$x$ is less or equal to $Kn$. \par 

 $\imath\imath$) $x$ is $(K+1)$-power free
(i.e. $u^{K+1}\in  \L(x)$ if and only if
$u=\emptyset$). \par 
 
$\imath\imath\imath$) 
For all $ u\in  \L(x)$ and for all $w \in\RR_{u}$ we have ${1\over
K}|u| < |w| $. \par  
                
 $\imath\nu$) For all $u\in
 \L(x)$, $\card ( \RR_u) \leq K(K+1)^2$.}

 \medskip

{\bf Proof.} We begin with a remark. Let $n$ be a
positive integer and $u\in  \L(x)$ a word of length $(K+1)n-1$. Let
$v\in \L(x)$ be a word of length $n$. The difference between two
successive occurrences of $v$ is less than $Kn$, consequently $u$
has at least one occurrence of $v$. We have proved that:                   
For each $n$, all words of length $n$ occurs in each word of length
$(K+1)n-1$. From this remark we deduce $\imath$).

Let $u\in  \L(x)$ be a word such that $u^{K+1}\in  \L(x)$.
Each factor of $x$ of length $|u|$ occurs in $u^{K+1}$. But in 
$u^{K+1}$ occurs at the most $|u|$ distinct factors of length $|u|$
of $x$. This contradicts the aperiodicity of $x$.

Assume there exist $u\in  \L(x)$ and $w\in  \RR_u$
such that $ |u|/K \geq |w|$. The word $w$ is a return word to $u$
therefore $u$ is a prefix of $wu$. We deduce that $w^{K}$ is a
prefix of $u$. Hence $w^{K+1}$ belongs to $ \L(x)$ because
$wu$ belongs to $ \L(x)$. Consequently $w=\emptyset$ and 
$\imath\imath\imath$) is proved.

Let $u$ be a factor of $x$ and $v \in \L(x)$ be a word of
length $(K+1)^2|u|$. Each word of length $(K+1)|u|$ occurs in $v$,
hence each return word to $u$ occurs in $v$. It follows from 
$\imath\imath\imath$) that in $v$ will occur at the most 
$K(K+1)^2|u|/|u|= K(K+1)^2$ return words to $u$, which proves
$\imath\nu$).  \cqfd

\medskip

We say that a subshift is {  \sl linearly recurrent} 
(resp.{\sl $K$-power free}) if it
is minimal and contains a linearly recurrent (resp. {\sl $K$-power
free}) sequence. If a subshift $(Y,T)$ is
linearly recurrent (resp. $K$-power free) then by minimality all
sequences belonging to $Y$ are linearly recurrent (resp. $K$-power
free). The following proposition was first proved in [Du1]. 

\medskip

{\bf Proposition 25.} { \sl All
substitution subshifts are linearly recurrent.} 

\medskip

{\bf Proof.}  Let $\tau $ be a substitution on $A$.
Let $u$ be a word of $ \L(\tau)$ and $v$ be a return word to $u$.
A well-known property of substitutions (see [Qu]) asserts that there
exists a constant $C$ such that for all positive integers $k$ 

$$ S_k= \sup \{ |\tau^k (a)| ; a\in A\}
 \leq C \inf \{ |\tau^k (a)| ; a\in A\} = CI_k\ . $$

 Let $k$ be the smallest integer such that $I_k \geq
|u|$. The choice of $k$ entails that there exists a word $ab\in
 \L(\tau)$ of length 2 such that $u$ occurs in $\tau^k (ab)$.
Let $R$ be the largest difference between two successive occurrences
of a word of length 2 of $ \L (\tau)$. We have 

$$ |v| \leq R S_k \leq RCI_k \leq RCS_1 I_{k-1} \leq RCS_1 |u|\ . $$ 

The subshift
spanned by $\tau$ is linearly recurrent with constant $RCS_1$.
 \cqfd                   

\medskip

{\bf Remarks.} Let $\varphi$ be a factor map from the subshift
$(X,T)$ on the alphabet $A$ onto the subshift $(Y,T)$ on the
alphabet $B$. A classical result (Theorem 6.2.9 in [LM]) asserts
that $\varphi$ is a {\sl sliding block code}. That is to
say there exists a non-negative integer $r$ and {\sl $r$-block map} $f: A^{2r+1}
\rightarrow B$ such that $(\varphi (x))_i = f(x_{[i-r,i+r]})$ for
all $i\in \Z$ and $x\in X$. We shall say that $f$ is a {\sl $r$-block map
associated to $\varphi$}. 

When $r=0$ (i.e. when $f$ is a map from $A$ to $B$), then $\varphi$ is the restriction to $X$ of the extension from $A^{\Z}$ to $B^{\Z}$ of $f$ introduced in subsection 2.1.

Suppose now that $r$ is a positive integer. If $u= u_0 u_1 \cdots u_n$ is a word of length $n\geq 2r+1$ we define
$\overline{f}(u)$ by $(\overline{f}(u))_i = f(u_{[i,i+2r]})$, 
$i\in \{ 0,1,\cdots n-2r\}$.
 Let $C$ denote the alphabet $A^{2r+1}$ and $Z=\{
((x_{[-r+i,r+i]}) ; i\in \Z)\in C^{\Z }; (x_n ; n\in  
\Z ) \in X \}$. It is easy to check that the subshift $(Z,T)$
is isomorphic to $(X,T)$ and that $\overline{f}$ induces a 0-block map from $C$
onto $B$ which defines a factor map from $(Z,T)$ onto $(Y,T)$.

\medskip

{\bf Proposition 26.} { \sl  Let $(X,T)$
and $(Y,T)$ be aperiodic subshifts such that $(Y,T)$ is a factor
of $(X,T)$. Then \par 
 
$\imath$) If $(X,T)$ is linearly recurrent with constant $K$
then $(Y,T)$ is linearly recurrent. Moreover there exists $n_0$ such
that: 
For all $u\in \L(Y) $, of length greater than
$n_0$, and for all $w\in  \RR_u$ we have $|u|/2K \leq |w|\leq
2K|u|$ and $\card (\RR_u) \leq 2K(2K+1)^2$. \par

$\imath\imath$) If $(Y,T)$ is $K$-power free, then $(X,T)$ is
$L$-power free for some $L$.} 

\medskip

{\bf Proof.} We denote
by $A$ the alphabet of $X$ and by $B$ the alphabet of $Y$. Let
$\varphi:(X,T) \rightarrow (Y,T)$ be a factor map. Let $f: A^{2r+1}
\rightarrow B$ be a block map associated to $\varphi$. 

Let $u$ be a word of $ \L(X)$ of length $|u|\geq 2r+1$ and $v$
be a word of $ \L(Y)$ defined by $f (u) = v$. We have 
$|v| = |u|-2r$. If $w$ is a return word to $v$ then 
$|w|\leq \max \{ |s|; s\in\RR_u \} \leq K|u| \leq K(|v| +2r)$.  The
subshift $(Y,T)$ is linearly recurrent with constant $K(2r+1)$.
Moreover: For all $v\in  \L(Y) $
such that $|v| \geq 2r$, and for all $w\in\RR_v$, $|w|\leq 2K|v|$.

Suppose there exists a word $\alpha$ of length $n\geq 2r$ such that
$\alpha^{2K+1}$ belongs to $ \L(Y)$. Each word of length $n$
occurs in $\alpha^{2K+1}$. But in $\alpha^{2K+1}$ occurs at the most
$|\alpha|=n$ words of length $n$. This contradicts the
aperiodicity of $(Y,T)$, therefore $\alpha$ is the empty word. The same
arguments we used to prove point $\imath\imath\imath$) of Theorem 24
complete the proof of $\imath$).

Assume that $(Y,T)$ is $K$-power free. Suppose that there exists a
word $u$, greater than $2r+1$, such that $u^{(K+1)}$ is a factor of $x$. Let $w$ be
$f(u u_{[0, 2r-1]})$. The word $w^K$ is a prefix of
$f(u^{(K+1)})$, hence $w$ is the empty word and the same for $u$. This
completes the proof. \cqfd                    

\medskip

{\bf Corollary 27.} { \sl Let $(X,T)$ and $(Y,T)$ be two isomorphic aperiodic subshifts. Then:

$\imath$) If $(X,T)$ is linearly recurrent, then $(Y,T)$ is linearly recurrent.

$\imath$$\imath$) If $(X,T)$ is $K$-power free for some $K$, then $(Y,T)$ is $K$-power free.
}

\medskip

{\bf 7.2. Chains of factors.}

\medskip

A {  \sl chain of factor maps of length $n$}
is a finite sequence $(\gamma_i : (Y_{i+1} , T_{i+1})\rightarrow (Y_{i} , T_{i}) ; 0\leq i\leq n-1)$ of factor maps of minimal Cantor systems :
$$ 
(Y_{n},T_n)\; \fleche{ \gamma_{n-1}} \;
(Y_{n-1},T_{n-1})\;  \fleche{ \gamma_{n-2}}\; 
\cdots\; \fleche{\gamma_1}(Y_1,T_1)\;
\fleche{\gamma_0} \; (Y_0,T_0)\ . 
$$
The notion of {\sl subchain} is implicitly defined. A chain is {\sl proper} when none of the $\gamma_i$'s is are isomorphisms.

Let $(\gamma_i : (Y_{i+1} , T)\rightarrow (Y_{i} , T) ; 0\leq i\leq n-1)$ be a chain of subshift factor maps. A remark in subsection 7.1 implies that by substituting, if needed, some $(Y_i,T)$ by an
isomorphic subshift, we can suppose that for each $i\in \{ 0,\cdots
, n-1\}$ the factor map $\gamma_i$ is
defined by a $0$-block map. 

\medskip

{\bf Proof of Theorem 3.} The
subshift $(Y,T)$ on the alphabet $A$ is linearly recurrent with
constant $K$. Let $(\gamma_i : (Y_{i+1} , T)\rightarrow (Y_{i} , T) ; 0\leq i\leq n-1)$ be a proper chain of aperiodic subshift factor maps of length $n$ such that $(Y_{n},T) = (Y,T)$. Substituting, if
needed, some $(Y_i,T)$, $0\leq i\leq n-2$, by an appropriate isomorphic subshift we
can assume that the factor maps $\gamma_{i}$, $0\leq i\leq n-2$, are given by $0$-block maps $\overline{\gamma}_i$ from the alphabet of $Y_{i+1}$ to the alphabet of $Y_i$. We do not change $(Y,T)$ and $\gamma_{n-1}$ because
it would change the constant $K$. It follows from Proposition 26
that there exists a positive integer $n_0$ such
that:                           
 For all $0\leq i \leq n-1$, for all
$v\in  \L(Y_i)$ such that $|v| \geq n_0$ and for all $w\in 
\RR_v$ we have $|v|/2K \leq |w| \leq 2K |v|$. We choose a
word $v_{n-1}$ of $ \L(Y_{n-1})$ such that $|v_{n-1}|\geq n_0$.

Let $ i \in \{0,\cdots ,  n-2 \}$. The factor map $\gamma_i \gamma_{i+1} \cdots
\gamma_{n- 2}$ is given by the 0-block-map $\psi =
\overline{\gamma}_i \overline{\gamma}_{i+1} \cdots
\overline{\gamma}_{n- 2}$. Let $v_i=\psi (v_{n-1})$. If $w$ is
a return word to $v_{n-1}$, then $\psi (w)$ is a concatenation of
return words to $v_i$. Therefore we can define a unique map $\lambda_i :                           
R_{v_{n-1}} \rightarrow R_{v_i}^*$ by $\phi_{v_i} \lambda_i = \psi
\phi_{v_{n-1}}$ (see 4.2). Let $y$ be a word of $ \L(Y_i)$. A word of length $n$ of
$ \L(Y_i)$ has at the most $2Kn/|y|$ occurrences of $y$. From
this we deduce that for all $a$ in $R_{v_{n-1}}$

 $$ |\lambda_i(a) | \leq 
{2K|\psi (\phi_{v_{n-1}} (a))|\over|v_i|} 
\leq  {2K|\phi_{v_{n-1}} (a)|\over |v_i|} \leq
{2K^2|v_{n-1}| \over |v_i|}  = 2K^2 . $$

Recall that $\card (R_{v_{i}}) \leq 2K(2K+1)^2$ (Theorem 24). 

The length of the image of each letter under the map $\lambda_i$, $0\leq i \leq n-2$, and the cardinality of the alphabets $R_{v_{n-1}}$ and $R_{v_i}$ are bounded independently of $n$ and $i$. Hence the set of maps $\lambda_i$, $0\leq i \leq n-2$, is bounded independently of $n$ and $i$, more precisely it is bounded by
$$
(2K(2K+1)^2)^{4K^3(2K+1)^2}=D. 
$$ 
Assume that $n$ is greater than
$D+2$, in this case there exists two integers $0\leq p < q\leq n-2$
such that $\lambda_p = \lambda_q$. Let $\rho$ be the $0$-block map  $\overline{\gamma}_p \overline{\gamma}_{p-1} \cdots
\overline{\gamma}_{q-1}$. It defines a factor map from $(Y_q,T)$
onto $(Y_p,T)$. We want to prove that $\rho ( \RR_{v_q}) =
 \RR_{v_p}$. 

The definition of $v_p$ and $v_q$ implies that $\rho (v_q) = v_p$.
Let $w$ be a return word to $v_q$. The word $\rho (w)$ is a
concatenation of return words to $v_p$. Suppose that $\rho (w)$ is a
concatenation of at least two return words to $v_p$. There exists
$b\in R_{v_{n-1}}$ such that $\overline{\gamma}_q\cdots
\overline{\gamma}_{n-2} (\phi_{v_{n-1}} (b)v_{n-1})$ has an occurrence of
$w v_q$. Consequently the number of occurrences of return words to
$v_p$ in $\overline{\gamma_p} \cdots \overline{\gamma}_{n-2}
(\phi_{v_{n-1}}(b)v_{n-1})$ is greater than the number occurrences of return
words to $v_q$ in $\overline{\gamma_q} \cdots
\overline{\gamma}_{n-2} (\phi_{v_{n-1}}(b)v_{n-1})$. That is to say,
$|\lambda_q (b)| < |\lambda_p (b)|$, which contradicts that
$\lambda_p = \lambda_q$. Therefore $\rho ( \RR_{v_q}) \subset
 \RR_{v_p}$.

Let $w$ be a return word to $v_p$. If there is no return word $x$ to
$v_q$ such that $\rho (x) = w$, then as above we can prove that
there exists a letter $b$ belonging to $R_{v_{n-1}}$ such that $|\lambda_q
(b)| < |\lambda_p (b)|$. We obtain $\rho ( \RR_{v_q}) = 
\RR_{v_p}$. 

Moreover $\card (  \RR_{v_q} ) = \card ( \RR_{v_p})$,
hence the letter to letter morphism $\rho$, restricted to the set of
concatenations of return words to $v_{q}$, is a bijection from the
set of concatenations of return words to $v_{q}$ onto the set of
concatenations of return words to $v_{p}$. This induces an
isomorphism from $(Y_q,T)$ onto $(Y_p,T)$ and proves that $\gamma_p\gamma_{p+1}\cdots\gamma_{q-1}$ is an isomorphism. This contradicts the fact that $n$ is the length of $(\gamma_i : (Y_{i+1} , T)\rightarrow (Y_{i} , T) ; 0\leq i\leq n-1)$.\cqfd          

\medskip

{\bf Proposition 28.} { \sl  Let $(Y,S)$
be a Cantor system which is a factor of a linearly
recurrent subshift. Then $(Y,S)$ is either isomorphic to a
subshift or to an odometer.} 
\medskip

{\bf Proof.}  Let $(X,T)$
be a linearly recurrent subshift with constant $K$ and $\lambda:
(X,T) \rightarrow (Y,S)$ be a factor map. Let $(\P_n = \{
P(n,1),\cdots ,P(n,p_n) \} ; n\in \N ) $ be a sequence of partitions
of $Y$ such that: \par 

 $\imath$) Each element of
these partitions is clopen.\par      

 $\imath\imath$)
$\P_n \prec \P_{n+1}$.\par                  
$\imath\imath\imath$) The
sequence $(\P _n ; n\in \N)$ separates the points of
$Y$. 

\medskip

 Let $n$ be a positive integer. We set $A_{n} = \{
1,\cdots ,p_n \}$ and we define a map $\delta: Y \rightarrow A_{n}$
by $\delta(x) = i$ if and only if  $x$ belongs to $P(n,i)$. Let
$\varphi_{n}:Y \rightarrow A_{n}^{\Z}$ be the map
defined by $(\varphi_{n} (x))_i = \delta (S^i(x))$ for all $i\in  
\Z $. We set $Y_n = \varphi_n(Y)$. The map $\varphi_n$ is a
factor map from $(Y,S)$ onto $(Y_n,T)$. 

Property $\imath\imath$) implies that there exists a $0$-block map from $A_{n+1}$ onto $A_n$,
which induces a factor map $\gamma_n:(Y_{n+1},T)
\rightarrow (Y_{n},T)$ satisfying $\gamma_n \varphi_{n+1} =
\varphi_n$. We have the following scheme:                           
$$ 
(X,T) \; \fleche{\lambda} \; (Y,S) \;
\fleche{\varphi_{n+1}} \;
(Y_{n+1},T)\; \fleche{ \gamma_n} \;
(Y_n,T)\;  \fleche{ \gamma_{n-1}}\; 
\cdots\; \fleche{\gamma_1}(Y_1,T)\;
\fleche{\gamma_0} \; (Y_0,T)\ . 
$$
Let $\Omega$ be the set $\{ (x_0,x_1,\cdots ) \in
\Pi_{i=0}^{+\infty} Y_i ; \gamma_{i}(x_{i+1}) = x_i , i\in \N                           \}$.
It is the inverse limit of the sequence $((Y_i,T),\varphi_i ; i\in
\N)$. We endow $\Omega$ with the topology
induced by the infinite product topology. Endowed with this topology
$\Omega $ is a Cantor set. Clearly, the map $U:\Omega\rightarrow
\Omega$ defined by $U((x_0,x_1,\cdots )) = (Tx_0,Tx_1,\cdots )$ is a
homeomorphism. The reader can check that $(\Omega , U)$ is a minimal
Cantor system. In what follows we prove that $(\Omega,U)$ is
isomorphic to $(Y,S)$.

Let $\gamma: (Y,S)\rightarrow (\Omega , U) $ be the map defined by
$\gamma (x) = (\varphi_0 (x),\varphi_1 (x),\cdots )$. It is easy to
see that $\gamma$ is a factor map. It remains to prove that $\gamma
$ is one-to-one. Let $x$ and $y$ be elements of $Y$. It follows from
property $\imath\imath\imath$) 
that there exists an integer $n$ and $1\leq i<j\leq p_n$ such that
$x\in P(n,i)$ and $y\in P(n,j)$. Hence we have $\varphi_n (x)_0 =
i\neq j = \varphi_n (y)_0$. This implies that $\gamma (x)\neq
\gamma (y)$. Therefore $\gamma$ is an isomorphism.

Taking the sequence of partitions $(\P_n ; n\geq l)$ for some
$l\in \N $, if needed, we can suppose that either \par                   
$\imath$) none of the subshifts $(Y_n,T)$ are periodic, or
\par                   
$\imath\imath$) each subshift $(Y_n,T)$ is periodic. \par 
          
Assume that
$(\P_n ; n\in \N)$ fulfills $\imath$). The subshift
$(X,T)$ is linearly recurrent; consequently Theorem 3 implies
that there exists $n_0\in \N  $ such that for all $n\geq n_0$
the subshift factor map $\gamma_n : (Y_{n+1},T)\rightarrow (Y_n,T)$ is an isomorphism. Then it is
easy to check that $(Y_{n_0},T)$ is isomorphic to $(\Omega ,U)$ and
{  \sl a fortiori} to $(Y,S)$.

Suppose $\imath\imath$). For each $n\in\N$ the subshift is
finite and we set $k_n = \card Y_n$. Clearly $(Y,S)$ is isomorphic to
the odometer with base $(k_n ; n\in \N)$.\cqfd  
    
\medskip

{\bf 7.3. Cantor factors of substitution subshifts.}

\medskip

In this part we give a condition for a subshift to
be isomorphic to a substitution subshift. Then we prove that under
some assumptions a Cantor factor of a substitution subshift is
isomorphic to a substitution subshift. We conclude this section with the proof
of Theorem 4.

\medskip

{\bf Proposition 29.} {\sl Let $x$ be a uniformly
recurrent sequence. Let $u=x_{[0,n]}$ and $v=x_{[0,l]}$, $n<l$, such
that each $wu$, $w\in  \RR_{u}$, occurs in each
return word to $v$. If $\D_{u} (x) =\D_{v} (x)$ then
\par                 
$\imath$) $\D_{u}(x)$ is a substitution fixed
point and \par
$\imath\imath$) the subshift spanned by $x$
is isomorphic to a substitution subshift.}

\medskip

{\bf Proof.} We have $\D_{u} (x) = \D_{v} (x)$ hence $R_{u} =
R_{v} = R$. It follows from the hypothesis that each return word to
$v$ is a concatenation of return words to $u$. This allows us to
define a morphism $\zeta: R \rightarrow
R^+$ by $\phi_{u} \zeta = \phi_{v}$. The reader can check it is a
substitution. Moreover, we have 
$$ \phi_{u} \zeta(\D_{u}(x)) =\phi_{v}(\D_{u}(x)) =
\phi_{v}(\D_{v}(x)) = x. $$ 
The unicity
of $\D_{u}(x)$ implies that $\zeta(\D_{u}(x)) =\D_{u}(x)$. 
This prove $\imath$).

We have $x = \phi_{u}(\D_{u}(x))$ and Lemma 17 asserts that $\phi_{u}(R)=\RR_{u}$ is a circular
code. Proposition 23 completes the proof. 
\cqfd                   

\medskip

{\bf Proof of Theorem 4. }  Let $(Y,S)$
be a minimal Cantor system which is the factor of a subshift spanned
by a substitution $\tau: A\rightarrow A^+$. Two cases appear
(Proposition 28):  \par 
  
 $\imath$) $(Y,S)$ is
isomorphic to a subshift $(W,T)$. \par

$\imath\imath$)  $(Y,S)$ is isomorphic to an odometer. 

\medskip

Assume that $(Y,S)$ fulfills $\imath$). Let $B$ be the
alphabet of $W$. Let $\varphi: (X_{\tau} , T) \rightarrow(W,T)$ 
be a factor map and $f: A^{2r+1}\rightarrow B$  an associated
block-map. Consider $A^{2r+1}$ as an alphabet that we denote $C$. Let $X$
be the set $\{ ((x_{[n-r , n+r]}) ; n\in \Z ) \in C^{\Z};
(x_n ; n\in\Z) \in X_{\tau} \}$.
The subshift $(X,T)$
is isomorphic to $(X_{\tau},T)$ and $f$ induces a $0$-block map from
$C$ to $B$ which defines a factor map
$\gamma: (X,T)\rightarrow (W,T)$. In [Qu]
is defined a substitution $\sigma:                           
C\rightarrow C^+$ in the following way  
$$ 
\sigma ((c_1 \cdots
c_{2r+1})) = (d_1 d_2 \cdots d_{2r+1}) ( d_2 d_3 \cdots d_{2r+2})
\cdots (d_{|\tau (c_1)|} d_{|\tau (c_1)|+1} \cdots d_{|\tau
(c_1)|+2r}), 
$$ 
where $\tau (c_1 \cdots  c_{2r+1}) = d_1 d_2 \cdots
d_{n}$. Let $(y_n ; n\in   \Z)$ be a fixed point
of $\tau$, and remark that $y^{'}=((y_{n} \cdots y_{n+2r});
n\in \Z )$ is a fixed point of $\sigma$
and belongs to $X$. Therefore $\sigma $ spans $X$. In [Du1] it is
proved that the image $z$ by a $0$-block map of a one-sided
substitution fixed point has a finite number of derivative
sequences. The proof can be adapted to bi-infinite sequences to
obtain that the set $\{ \D_{u} (z)\;;\;
u=z_{[0,n]}, n\in \N\}$ is finite. We
apply this result to the uniformly recurrent sequence $ \gamma
(y^{'})$ which spans $W$: There exists
$u$ and $v$ fulfilling the hypothesis of Proposition 29. Hence $(W,T)$, and { 
\sl a fortiori} $(Y,S)$, is isomorphic to a substitution
subshift. 

\medskip

Now we suppose that $(Y,S)$ is isomorphic to the odometer with base
$(p_n;n\in \N)$. To prove that this odometer has
a stationary base it suffices to prove that the subset $\Gamma$ of
prime numbers dividing some $p_n$ is finite.

We choose $r$, $l\in A$ satisfying the requirement $\imath$),
$\imath\imath$)
and $\imath\imath\imath$) in subsection 2.3.2. Let $u$ be a return
word to $r.l$. Let $p$ be a prime number of $\Gamma$. The system
 $( \Z /p\Z, x\mapsto x+1\; {\rm (mod} \; p{\rm )})$ 
is a factor of $(X_{\tau}, T)$. Consequently there exists a positive
integer $n_0$ such that $p$ divides $|\tau^n (u)|= \sum_{a\in
A}{(M(\tau)^n(|u|_a ; a\in A) )_a}$ for all $n\geq n_0$. Let $P(X) =
\sum_{i=0}^{|A|} c_i X^i$ be the characteristic polynomial of
$M(\tau)^n$. We know that $c_0 = {\rm det}( M(\tau)^n)$ and, by the
Theorem of Cayley-Hamilton, that $P(M(\tau)^n)(|u|_a ; a\in A)=0$. It follows that
$p$ divides $c_0 |u|$, hence $p$ divides $|u|{\rm det} M(\tau)$. We conclude that
the set $\Gamma$ is finite. 

\medskip

Suppose now that $(Y,S)$ is a minimal Cantor system which is a factor of an odometer with stationary base $(p,q,q,\cdots)$. Let $(\P_n ; n\in \N)$ be a nested sequence of partitions of $Y$. 
Let $n\in \N$. Let $(Y_n,T)$ be the subshift associated to the partition $\P_n$ (see the proof of Proposition 28). It is easy to check that this subshift is finite and its period $k_n$ divides $pq^i$ for some $i$. Hence the set of prime numbers dividing some $k_n$ is finite and $(Y,S)$ is isomorphic to an odometer with stationary base.
\cqfd

\bigskip

%% file: section8.tex
{\gros 8. One--sided systems.}

\medskip

The motivation of this section lies in the fact that people working with
substitutions often consider only one--sided fixed points, bypassing the
difficulty of admissibility which  occurs when considering two--sided
ones. One is lead to consider one--sided substitution dynamical
systems, defined as the closure of the orbit of some one--sided fixed point
in $A^\N$. After the basic definitions, we shall introduce here the
dimension group of an one--sided system, and explain its relations with
the dimension group of  (ordinary, i.e. two-sided) dynamical systems. It will
lead us to an alternative method --- which could be called one--sided
--- to compute the dimension group of a substitution dynamical system;
more precisely, it uses the return words on $\ell$ only,  instead of the return
words on $r.\ell$ used in Section 5.

\medskip

{\bf 8.1. Definitions}

\medskip

{\bf Definition 14.} {\sl A {\rm one--sided dynamical system} $(X,T)$ is a
compact metric space $X$ endowed with a map $T\colon X\to X$ which is continuous and onto.}

\medskip

Recall that we say that $(X,T)$ is a dynamical system when in addition $T$ is also one-to-one. The most classical examples of one--sided systems are
one--sided subshifts:

\medskip

{\bf Definition 15.} {\sl Let $A$ be an alphabet, $T$ the shift on $A^\N$
defined by $(Tx)_n=x_{n+1}$ for all $x\in A^\N$ and all $n\in \N$, and $X$
a closed subset of $A^\N$ with $TX=X$. We say that $(X,T)$ is a {\rm one--sided
subshift on} $A$.}

\medskip

The one--sided subshift $(X,T)$ is said to be {\sl minimal} if the
only closed sets $F$ of $X$ with $TF=F$ are $\emptyset$ and $X$; it implies
that the only closed sets $F$ with $TF\subset F$ are $\emptyset$ and $X$.

\medskip

Some of the one-sided dynamical systems can also be associated to ordered
Bratteli diagrams: A simple ordered Bratteli diagram  is said to be {\sl
semi--proper} if it has only one minimal path. For such a diagram
$\B$ and $x\in\xb$, we define $\vb x$ as usual if $x$ is not maximal, and
$\vb x$ to be the minimal path $\xmin$ if $x$ is maximal. The system
$(\xb,\vb)$ is then a minimal one--sided system. But such a system is very
particular: it has only one point --- namely $\xmin$ --- with several
preimages. In fact, there is no clear method to associate a Bratteli diagram
to minimal one--sided Cantor systems in general. 

\medskip

{\bf 8.2. The natural extension.} 

\medskip

Given a one--sided system $(X,T)$ there
exists a dynamical system $(\tilde X,\tilde T)$ and a factor map $\pi
:\tilde X\to X$ such that:

\medskip

{\leftskip 1cm 
For every dynamical system $(Y,S)$ and every factor map
 $\phi\colon Y\to X$ there exists a unique factor map $\tilde\phi\colon
Y\to \tilde X$ with $\pi\circ\tilde\phi=\phi$.\par}

\medskip

 The triple $(\tilde X,\tilde T,\pi)$ with this property is unique up to
isomorphism in an obvious sense. It is called {\sl the natural extension of}
$(X,T)$. 

Let us
recall one of the possible constructions. Let $X^\Z$ be endowed with the
product topology, and $\tilde X$ be the subset of $X^\Z$ consisting of the
points $x=(x_n\semi n\in \Z)$ with $x_{n+1}=Tx_n$ for all $n\in\Z$.
$\tilde X$ is closed in $X^\Z$, it is not empty because $T$ is onto, and it is
invariant under the shift of $X^\Z$. Let $\tilde T$ be the restriction of this
shift to $\tilde X$, and $\pi\colon \tilde X\to X$ be defined by $\pi(x)=x_0$.
$(\tilde X,\tilde T,\pi)$ satisfies all the announced properties, thus can be
called the natural extension of $(X,T)$. From this construction it follows
that:

{ \leftskip 2cm \sl
$\imath$)  $\tilde X$ is a Cantor set if $X$ is a Cantor set.

$\imath\imath $) $(\tilde X,\tilde T)$ is minimal if and only
if $(X,T)$ is minimal.

$\imath\imath\imath$) For every clopen set $U$ of $\tilde X$
there exists $n\geq 0$ such that $\pi(\tilde{T}^{-n}U)$ is clopen and 
$\tilde{T}^{-n}U=\pi^{-1}\pi(\tilde{T}^{-n}U)$.
\par
}

\medskip

An alternative construction is possible when $(X,T)$ is a one--sided
subshift on some alphabet $A$: Let $\L(X)$ be the language of $X$, i.e. the
set of words on $A$ which are factors of $x$ for some $x\in X$; let
 $(\tilde X,\tilde T)$ be the (two--sided) subshift associated to this
language: $\tilde X$ is the set of $y\in A^\Z$ every factor of which belongs
to $\L(X)$, and $\tilde T$ be the restriction of the shift of $A^\Z$ to
$\tilde X$. The projection $\pi\colon A^\Z\to A^\N$ maps $\tilde X$ onto
$X$, and $(\tilde X,\tilde T,\pi)$ is the natural extension of $(X,T)$.

\medskip

{\bf 8.3. The dimension group of a one--sided system.}

\medskip

In this subsection we assume the spaces we encounter are Cantor sets.

A dimension group can be associated to every Cantor minimal one--sided
system $(X,T)$, exactly as for dynamical systems [HPS]: Let
$\partial_T C(X,\Z)$ be the group of coboundaries, i.e. the set of
functions which can be written $g-g\circ T$ for some $g\in C(X,\Z)$.
$K^0(X,T)$ is the quotient group $ C(X,\Z)/\partial_T C(X,\Z)$;
$K^0_+(X,T)$ is the image of $ C(X,\Z_+)$ in this quotient; and
$\one$ the image of the constant function $1$. $(K^0(X,T),K^0_+(X,T))$ is
an ordered group, and $\one$ an order unit. The triple  
$(K^0(X,T),K^0_+(X,T),\one)$ is called the {\sl dimension group of the
one--sided system} $(X,T)$. (That it actually is a dimension group follows from Proposition 32.)

The dimension groups of one--sided systems share some of the properties
of those of dynamical systems. In particular, as in [GW], we have:

\medskip
{\bf Lemma 30.} {\sl Let $\phi\colon (X,T)\to (Y,S)$ be a factor map between
two one--sided systems. Then the corresponding homomorphism
$\phi^*\colon K^0(Y,S)\to K^0(X,T)$ is one-to-one,  and 
$ \phi^*(K^0_+(Y,S)) = K^0_+(X,T)\cap \phi^*(K^0(Y,S))$.
}

\medskip

Moreover, we have:

\medskip

{\bf Proposition 31.} {\sl Let $(\tilde X,\tilde T,\pi)$ be the natural extension
of the one--sided system $(X,T)$. Then $\pi^*\colon K^0(X,T)\to
K^0(\tilde X,\tilde T)$ is an isomorphism of ordered groups with order
units.}

\medskip

{\bf Proof.} By Lemma 30 above, we have only to prove that $\pi^*$
is onto.  Let $\alpha\in K^0(\tilde X,\tilde T)$, and $f\in C(\tilde
X,\Z)$ a representative of $\alpha$. From the remark $\imath\imath\imath$) of subsection 8.2 it follows
that there exists an integer $n\geq 0$ such that $f\circ T^n$ factorizes through
$X$, i.e. is equal to $g\circ\pi$ for some  $g\in C(X,\Z)$.
 If $\beta$
is the class of $g\in K^0(X,\Z)$, then $\pi^*(\beta)$ is the class of $f\circ
T^n$ in $K^0(\tilde X,\tilde T)$, which is $\alpha$. \cqfd

\medskip

Let $(Y,S)$ be a Cantor minimal one--sided system, $(X,T,\pi)$ its
natural extension. When we want to compute the dimension group of
$(Y,S)$ or of $(X,T)$,  we are often lead to use  a
sequence of Kakutani--Rohlin partitions of $X$ which doesn't satisfy
neither the hypotheses $\imath\imath\imath$) nor $\imath\nu$) of subsection 1.3. We have to work with weaker
conditions, and the next result can't be obtained  by simple modifications
of the proofs of [HPS]. As the use of Bratteli diagrams seems of little help,
we prefer use the formalism of [GW], associating directly a
dimension group to the sequence of partitions 1.4.

\medskip

{\bf Theorem 32.} {\sl Let $(X,T,\pi)$ be the natural extension of the
one--sided Cantor minimal system $(Y,S)$.  Let $(\P_n = \{ T^j B_{n,k} ; k\in A_n  , 0\leq j <h_{n,k} \} ; n\in \N)$ be a nested
sequence of clopen Kakutani--Rohlin partitions of $X$, where $B_n$ is the base of $\P_n$, such that:\par

{\leftskip 1cm
 $\imath$) If $x,x'\in X$ are such that $\pi(x)\neq\pi(x')$, then for
$n$ large enough\par they belong to different elements of $\P_n$.\par 
 $\imath\imath$) For each $n$ there exists $p>n$ and $k\in A_n$ such that 
$B_p\subset B_{n,k}$\par
}

 Then $K^0(X,T)$ is  isomorphic as ordered group with order unit to the
dimension group   $K(\P_n\semi n\geq 1)$ defined in subsection 1.4. }

\medskip
{\bf Proof.} We use here freely the notations of subsection 1.4.

 Let $x,x'$
be two points in $\bigcap_p B_p$, and suppose that $\pi(x)\neq\pi(x')$. By 
$\imath$) there exists $n$ such that $x$ and $x'$ belong to different
elements of $\P_n$.  Let $p$ and $k$ be as in $\imath\imath$); $x$ and $x'$
both belong to $B_{n,k}$, and we get a contradiction. We have proved:

{\hfill $\displaystyle\pi\bigl(\bigcap_pB_p\bigr)$ consists of one
point.\hfill(1)}

We claim that
$$
\inf_{k\in A_n}\;h(n,k)\to +\infty\hbox{ when }n\to +\infty
\ .\eqno(2)
$$
Suppose it is false. There exists an integer $J>0$ and, for infinitely many
values of $n$, a point $x_n\in B_n$ and an integer $j_n$ with $0<j_n \leq J$
and $T^{j_n}x_n\in B_n$. By compactness, there exists a point $x\in\bigcap_n
B_n$ and an integer $j>0$ such that $T^jx\in\bigcap_nB_n$. By (1), 
$S^j\pi(x)=\pi(T^jx)=\pi(x)$, and it is impossible because $(Y,S)$ is minimal
and $Y$ is infinite: our claim is proved.

We prove now  that $\lambda$  maps $K^+(\P)$ onto $K^0_+(X,T)$.
Let $f\in C(X,\Z^+)$.
By Proposition 31, $f$ is cohomologuous to $g\circ\pi$ for some 
$g\in C(Y,\Z^+)$.
 From $\imath$) it follows by compactness that, for every clopen set
 $A$  of $Y$, $\pi^{-1}(A)$ is a union of elements of $\P_n$
for $n$ large enough. 
Thus  $g\circ \pi\in C_n^+$ for $n$ large enough. It follows that
$\lambda(K_n^+)=K^0_+(X,T)$, and it implies  immediately  that $\lambda$
is onto.

It remains to prove that $\lambda$ is one-to-one.
Let $\alpha$ be an element of $K(\P)$ with $\lambda(\alpha)=0$.
$\alpha$ is the image in $K(\P)$ of an element $\beta\in K_r$ for
some $r\geq 1$; let $f\in C_r$ be a representative of $\beta$. As
$\lambda(\alpha)=0$,  the definition of $\lambda$ implies that $f$ is the
coboundary of some   $g\in C(X,\Z)$.

 As in the proof of Proposition 31, there exist $i\geq 0$ and $h\in C(Y,\Z)$ such that 
$g\circ T^i=h\circ\pi$.

By (2), there exists $n\geq r$ such that $h(n,k)>i$ for every 
$k\in\ A_n$. 
Let $p$ and $k$ be  associated to $n$ as in $\imath\imath$).
 For $0\leq j<i$, $T^jB_p\subset T^j B(n,k)$, and $f$ is
constant on this set because it belongs to $\P_n$ and $f\in C_n$:
 $f$ is constant on $T^jB_p$ and $f\circ T^j$ is constant on $B_p$.

By (1),  there exists $q\geq p$ such that the
function $h$ is constant on $\pi(B_q)$, and it follows that $g$ is constant
on $T^iB_q$, and that $g\circ T^i$ is constant on $B_q$.  

As $f=g-g\circ T$,  
$$g = g\circ T^i+\sum_{j=0}^{i-1}f\circ T^j\; ;$$
 each of the functions on the right side of this expression is constant on
$B_q$, thus $g$ is constant on $B_q$. 

Moreover, $f$
belongs to $ C_q$ and,  for all   $k\in A_n$ and  all $x\in B_q$,

$$\sum_{j=0}^{h(q,k)-1} f(T^jx)=g(x)-g(T^{h(q,k)}x)=0$$
because $x$ and $T^{h(q,k)}x$ both belong to $B_q$ and $g$ is constant on
this set. Thus $f$ belongs to  $H_q$, and the projection of $f$ in $K_q$ is
$0$. But this projection is equal to $i_{q-1}\circ\cdots\circ i_r(\beta)$.
From the definition of a direct limit, it follows that $\alpha=0$.\cqfd

\medskip
{\bf An Application: one--sided return words.}
In section 4.3. we have considered a point $x$ in  a minimal subshift
$(X,T)$, and we have seen how to compute the dimension group of $(X,T)$
by using the return words on $x_{[-n,-1]}\,.\,x_{[0,n)}$.
We could make the same construction but with the return words
on $x_{[0,n)}$, and it would produce a sequence $(\QQ'_n)$ of partitions
instead of $(\QQ_n)$.
 Let $Y$ be the projection of $X$ in $A^\N$ and $S$ the shift on $A^\N$:
$(X,T)$ is the natural extension of $(Y,S)$.  
The sequence $( \QQ'_n)$
satisfy the hypotheses of Theorem 32, thus it can be used to compute
the dimension group of the system.

\medskip

{\bf 8.4. Application to substitution dynamical systems.}

\medskip

{\bf Definition 16.} {\sl We say that the substitution $\sigma$ on the
alphabet $A$ is {\rm left--proper} if there exist a letter $\ell\in A$ and an
integer $p>0$ such that $\ell$ is the first letter of $\sigma^p(a)$ for all
$a\in A$.}

\medskip

{\bf Corollary 33.} (Compare with Theorem 22 in 5.4) {\sl Let $\sigma$ be a
left--proper (primitive, aperiodic) substitution  on the alphabet $A$.
Then the dimension group of $(\xs,\ts)$ is isomorphic to the dimension
group of its matrix (as ordered groups with order units).}

\medskip

{\bf Proof.} Let $\pi$ be the natural projection of $A^\Z$ to $A^\N$, $S$ the one--sided shift on $A^\N$, and $Y=\pi(\xs)$; $(Y,S)$ is the
one--sided substitution dynamical system associated to $\sigma$, and 
$(\xs,\ts,\pi)$ is the natural extension of $(Y,S)$. We have only to check
that the sequence of partitions $(\P_n)$ introduced in subsection 2.4 satisfies
the two conditions of Theorem 32. But it is immediate, using the
same method as in the proof of the point $\imath\imath$) of Corollary 12. \cqfd

\medskip
The preceeding remarks lead to an alternative method to compute the
dimension group of an arbitrary substitution dynamical system. Let
$\sigma$ be a (primitive, aperiodic) substitution on the alphabet $A$, $z$
one of its one--sided fixed points, and $\ell=z_0$. There exists an
admissible two--sided fixed point $x$ with $x_{[0,+\infty)}=z$.
 We use the notations of
Section 4, following the same method as in subsection 5.2, except that we consider the
return words on $ \ell$ instead of $r.\ell$. Let $y=\D_{\ell}(x)$. As in subsection 5.2,
there exists a substitution $\rho$ on $R_{\ell}$ such that
$\phi_{\ell}\circ \rho=\sigma\circ \phi_{\ell}$, and $y$ is a fixed point of
$\rho$. $\rho$ is primitive and left-proper, and $(X_\rho,T_\rho)$ is
isomorphic to the system induced by $(\xs,\ts)$ on the cylinder set
$[\ell]$. From the Corollary 33 we get:

\medskip

{\bf Proposition 34.} {\sl With the definitions above, the dimension group of
$(\xs,\ts)$ is the dimension group of the matrix of $\rho$, except that the
order unit is the image of the vector $\bigl(|\phi_{\ell}(j)|\semi j\in
R_{\ell}\bigr)$.}

\medskip

{\bf Remark.} This method to compute the dimension group could be
called one-sided: The two-sided fixed point $x$ appears here only because
we have no defined return words for one-sided sequences.

%% file: biblio.tex
{\gros References:}

\indent

\bigskip

\item{[AP]} J.E. Anderson,I.F. Putnam, Topological invariants for substitution tilings
and their associated $C^*$-algebras, preprint.

\item{[Be]} J. Bellissard, Gap labelling theorems for Schr\"odinger's operators, Numb. Theo. and Physics, J. M. Luck, P. Moussa and M. Waldschmidt (Eds.), Springer proceedings in Physics 47 (1990).

\item{[BBG]} J. Bellissard, A. Bovier and J.-M. Ghez, Gap labelling theorems for one dimensional discrete Schr\"odinger operators, Rev. in Math. Phy. 1 (1992), 1-37.

\item{[BT]} E. Bombieri and J. E. Taylor, Quasicrystals, tilings and algebraic number theory: some preliminary connections, Cont. Math. 64 (1987), 241-264.

\item{[Br]} O. Bratteli, Inductive limits of finite-dimensional $C^*$--algebras, Trans. Amer. Math. Soc. 171 (1972), 195--234.

\item{[BHMV]}  
V. Bruy\`ere, G. Hansel, C. Michaux, R. Villemaire, Logic and $p$-recognizable sets of integers, Bull. Belgian Math. Soc. Simon Stevin vol. 1 (1994) 191-238.

\item{[CKMR]} G. Christol, T. Kamae, M. Mendes-France and G.
Rauzy, Suites Alg\'ebriques et Substitutions, Bull. Soc. Math. France 108
(1980), 401--419. 

\item{[Co1]} A. Cobham,  On the base-dependence of sets of numbers recognizable by finite automata, Math. Syst. Theo. 3 (1969), 186--192. 

\item{[Co2]} A. Cobham, Uniform tag sequences, Math. Syst. Theo. 6 (1972), 164--192. 

\item{[De1]} F. M. Dekking, Transcendance du nombre de Thue-Morse, C. R. Acad. Sci. Paris 285 s\'erie A (1977), 157--160.

\item{[De2]} F. M. Dekking, The spectrum of dynamical systems arising from substitutions of constant length, Zeit. Wahr., 41 (1978), 221--239.

\item{[De3]} F. M. Dekking, Recurrent sets, Adv. in Math. 44 (1982), 78--104.

\item{[Du1]} F. Durand, A characterization of substitutive sequences using return words, to appear in Discrete Mathematics.

\item{[Du2]} F. Durand, A generalization of Cobham's theorem, accepted by Math. Syst. Theo..

\item{[El]} G. A. Elliott, On the classification of inductive limits of sequences of semi-simple finite dimensional algebras, J. Algebra 38 (1976), 29--44.

\item{[Ef]} E. G. Effros, Dimensions and $C^*$--algebras, Conf. Board Math. Soc. Providence, R.I. (1981).

\item{[EHS]} E. Effros, D. Handelman and C.-L. Shen, Dimension groups and their affine representations, Amer. J. Math. 102 (1980), 385--407.

\item{[FM]} 
S. Ferenczi and C. Mauduit, Transcendency of numbers with low complexity expansion, preprint.

\item{[FMN]} S. Ferenczi, C. Mauduit  and A. Nogueira, Substitution dynamical systems: Algebraic characterization of eigenvalues, Ann. Scient. Ec. Norm. Sup. 29 (1996), 519-533.

\item{[Fo]} A. H. Forrest, $K$--groups associated with substitution minimal systems, to appear in Israel Journal.

\item{[GPS]} T. Giordano, I. Putnam and C. F. Skau, Topological orbit equivalence and $C^*$--crossed products, J. reine angew. Math. 469 (1995), 51--111.

\item{[GW]} E. Glassner and B. Weiss, Weak orbit equivalence of Cantor minimal systems, Internat. J. Math. 6 (1995), 559--579.

\item{[Go]}
W. H. Gottschalk, Substitution minimal sets, Trans. Amer. Math. Soc. 109 (1963), 467-491.

\item{[HL]} T. Harju and M. Linna, On the periodicity of morphisms on free monoids, Informatique th\'eorique et applications~/ Theoretical informatics and applications, vol. 20 n. 1 (1986), 47--54.

\item{[HM1]} 
G. A. Hedlund and M. Morse, Symbolic dynamics, Amer. J. Math. 60 (1938), 815-866.

\item{[HM2]} 
G. A. Hedlund and M. Morse,  Symbolic dynamics II. Sturmian trajectories, Amer. J. Math. 62 (1940), 1-42.

\item{[HPS]} R. H. Herman, I. F. Putnam and C. F. Skau, Ordered Bratteli diagrams, dimension groups and topological dynamics, Internat. J. Math. 3 (1992), 827--864.

\item{[Ho]} B. Host, Valeurs propres des syst\`emes dynamiques d\'efinis par des
substitutions de longueur variables, Ergod. Th. \& Dynam. Sys. 6 (1986), 529--540.

\item{[IK]} S. Ito and M. Kimura, On Rauzy fractal, Japan J. Indust. Appl. Math. 8 (1991), 461--486.

\item{[LM]} D. Lind and B. Marcus, An
Introduction to Symbolic Dynamic and Coding, Cambridge University
Press (1995).

\item{[LV]} A. N. Livshitz and A. M. Vershik, Adic models of ergodic theory, spectral theory, substitutions, and related topics, Adv. in Sov. Math. 9 (1992), 185--204.

\item{[LP]} J. H. Loxton and A. J. Van der Poorten,  Arithmetic properties of the solutions of a class of functional equations, J. reine angew. Math. 330 (1982), 159-172.

\item{[Ma]} J. C. Martin, Substitution minimal flows, Amer. J. Math. 93 (1971), 503--526.

\item{[MS]} F. Mignosi and P. S\'e\'ebold, If a D0L--language is $k$--power free then it
is circular, ICALP 1993, Lect. Notes in Comp. Sci. 700.

\item{[Mo1]} B. Moss\'e, Puissances de mots et reconnaissabilit\'e des points fixes
d'une substitution, Theo. Comp. Sci. 99 (1992), 327--334.

\item{[Mo2]} B. Moss\'e, Reconnaissabilit\'e des substitutions et complexit\'e des
suites automatiques, Bull. Soc. Math. France 124 (1996), 101--118.
 
\item{[Pa]} J.J. Pansiot, Decidability of periodicity for infinite words, Informatique
th\'eorique et applications~/ Theoretical informatics and applications, vol.
20 n. 1 (1986), 43--46.

\item{[Pe]} R. Penrose, The role of aesthetics in pure and applied mathematical research, Bull. Inst. Math. Appl. 10 (1974), 266--271.

\item{[Qu]} M. Queff\'elec, Substitution Dynamical systems, Lecture Notes in Math. 1294 (1987).

\item{[RS]} G. Rozenberg and A. Salomaa, The mathematical theory of {\rm L} systems, Academic Press (1980).

\item{[SBGC]} D. Shechtman, D. Blech, I. Gratias and J. V. Cahn, Metallic phase with long-range orientational order and no translational symmetry, Phys. Rev. Lett. 53 (1984), 1951-1953.

\item{[So]} B. Solomyak, On the spectral theory of adic transformations, Adv. Sov. Mat. 9 (1992), 217-230.

\item{[Ve]} A. M. Vershik, A theorem on the Markov periodical approximation in ergodic theory, J. Sov. Math. 28 (1985), 667--674.